\documentclass[10pt]{amsart}
\usepackage{latexsym}   
\usepackage{amssymb}    
\usepackage{amsmath}    
\usepackage{amsthm}     
\usepackage[all]{xy}
\usepackage{hyperref}
\newcommand{\pa}{\partial}\newcommand{\al}{\alpha}
\newcommand{\be}{\beta}

\newcommand{\Ga}{\Gamma}\newcommand{\del}{\delta}
\newcommand{\ep}{\epsilon}

\newcommand{\la}{\lambda}\newcommand{\La}{\Lambda}\newcommand{\om}{\omega}
\newcommand{\Om}{\Omega}
\newcommand{\si}{\sigma}\newcommand{\ti}{\tilde}

\renewcommand{\thefootnote}

\title[ The B\"{a}cklund transforms of Peterson's deformations of quadrics]
{The B\"{a}cklund transforms of Peterson's deformations of
quadrics}

\author[  Ion I. Dinc\u{a}]{Ion I. Dinc\u{a}}
\address{Faculty of Mathematics and Informatics,
University of Bucharest,  14 Academiei Str., 010014, Bucharest,
Romania}
 \email{dinca@gta.math.unibuc.ro}
\thanks{Supported by the University of Bucharest}

\begin{document}


\begin{abstract}
In trying to provide explicit deformations of quadrics the
starting point of our investigation is to use Bianchi's link
between real deformations of totally real regions of real
paraboloids and various totally real forms of the sine-Gordon
equation coupled with Bianchi's simple observation that the vacuum
soliton of these totally real forms of the sine-Gordon equation
provides precisely Peterson's deformations of such quadrics in
order to derive explicit B\"{a}cklund transforms of Peterson's
deformations of quadrics. Based also on Bianchi's approach of the
B\"{a}cklund transformation for quadrics via common conjugate
systems and in analogy to the solitons of the sine-Gordon equation
corresponding at the level of the geometric picture to the
solitons of the pseudo-sphere we propose a model for the solitons
of quadrics.
\end{abstract}

\maketitle

\tableofcontents \pagenumbering{arabic}

\section{Introduction}
During the years 1899-1906 the theory of deformation (through
bending) of general quadrics got the attention of geometers
(mainly Bianchi, Calapso, Darboux, Guichard, Peterson and
\c{T}i\c{t}eica); as a consequence of their results the classical
differential geometry of surfaces underwent a fundamental change.
This theory culminated with Bianchi's discovery in 1906 of the
{\it B\"{a}cklund} (B) transformation for general quadrics and the
{\it applicability correspondence provided by the Ivory affinity}
(ACPIA). However no explicit examples of deformations built on
Bianchi's approach exist in literature except mainly for the
solitons of the (pseudo-)sphere. As any other integrable system
one method to produce explicit solutions is to begin with the
vacuum soliton as seed and build its B transforms. However in our
case another seeds (namely Peterson's deformations of quadrics)
will be amenable to explicit computations of their B transforms.
The starting point of our investigation is to use Bianchi's link
from (\cite{B1}, ch VI) between real deformations of totally real
regions of real paraboloids and various totally real forms of the
sine-Gordon equation coupled with Bianchi's simple observation
that the vacuum soliton of these totally real forms of the
sine-Gordon equation provides precisely Peterson's deformations of
such quadrics to derive explicit deformations of quadrics.

The condition that a {\it conjugate system} (that is the second
fundamental form is missing mixed terms) on a quadric is common to
a Peterson's $1$-dimensional family of deformations of the quadric
is projective invariant; also the condition that the lines of
coordinates are planar is a projective invariant. On the complex
unit sphere such conjugate systems with planar lines of
coordinates are given by orthogonal systems of circles, that is
the axes of the two pencils of planes containing the circles are
polar reciprocal with respect to the sphere; according to Bianchi
this condition is projective invariant, so it is valid for all
quadrics. For general quadrics when one of the axes is a principal
axis for the quadric one can derive explicit formulae for
Peterson's $1$-dimensional family of deformations of quadrics.

The starting point of Bianchi's investigation was results of
Calapso, Darboux and Servant according to which for any real
deformation of a totally real region of a quadric the conjugate
system common to the deformation and quadric (any two surfaces in
a point-wise correspondence admit a common conjugate system) is
{\it isothermal-conjugate} system of coordinates on the
deformation (that is the second fundamental form is a multiple of
the identity, modulo some signs as required by curvature and
signature of the ambient space considerations): he introduced as
an auxiliary variable a (hyperbolic) angle. According to Bianchi
the (hyperbolic) angle between one of Peterson's conjugate system
lines and one of the lines of the isothermal-conjugate system is a
solution of some totally real form of the sine-Gordon equation for
general real paraboloids: this is the geometric link between the
sine-Gordon equation and totally real deformations of totally real
regions of real paraboloids.

Note also that Calapso in \cite{Ca} has completed Bianchi's
approach of the B transformation of deformations of
$2$-dimensional quadrics via common conjugate systems from
paraboloids to {\it quadrics with center} (QC), but his approach
for QC is different from Bianchi's outline. The condition that the
conjugate system on a quadric is a conjugate system on one of its
deformations was known to Calapso for a decade, but the
B\"{a}cklund transformation for general quadrics via the Ivory
affinity eluded Calapso since the common conjugate system is
a-priori best suited for the B transformation only at the analytic
level (which makes it also the best suited tool to provide
explicit examples).

In what concerns totally real deformations of totally real regions
with positive linear element of real paraboloids our main result
is to put Bianchi's machinery to work to churn out explicit
formulae for the B transforms of Peterson's deformations of such
quadrics up to including the third iteration of the B
transformation.

In what concerns totally real deformations of totally real regions
of other quadrics our main result is to complete Bianchi's elegant
approach of the B transformation via common conjugate systems and
then use this to churn out explicit formulae for the B transforms
of Peterson's deformations of such quadrics up to including the
third iteration of the B transformation. Calapso's approach from
\cite{Ca} (another completion of Bianchi's approach of the B
transformation via common conjugate systems to general quadrics)
is similar in the main ideas to Bianchi's approach but different
in the fact that he uses only the common conjugate system, without
paying attention to the change from the initial conjugate system
common to a Peterson's $1$-dimensional family of deformations of
quadrics. The totally real forms of the sine-Gordon equation are
replaced for quadrics with center by another equation

Once a case of a general QC (the general case) being solved, all
other complex types of quadrics should be amenable to explicit
computations of the B transforms of Peterson's deformations of
such quadrics by similar computations. Since multiplication by $i$
exchanges both the signature of the totally real surface and of
the ambient Lorentz space, from a totally real point of view one
needs only discuss deformations in
$\mathbb{R}^2\times\sqrt{\ep}\mathbb{R},\ \ep=\pm 1$ of quadrics
with positive linear element (there are for example isotropic
quadrics without center that cannot be realized as real quadrics,
but admit real deformations) and deformations in
$\mathbb{R}^2\times i\mathbb{R}$ of quadrics with linear element
of signature $(1,1)$.

The deformation problem for positive definite linear element is
elliptic for real deformations of surfaces of positive Gau\ss\
curvature and for totally real deformations in Lorentz spaces of
signature $(2,1)$ of surfaces of negative Gau\ss\ curvature and
hyperbolic otherwise.

Thus for the hyperbolic sine-Gordon and sinh-Gordon equation the
seed and the leaf will admit asymptotic lines and will be
applicable to the same totally region of the real quadric; for the
elliptic sine-Gordon and sinh-Gordon equation the seed and the
leaf  will not admit asymptotic lines and will be applicable to
different totally regions of the real quadric (the applicability
becomes ideal in Peterson's denomination); one needs composition
of B transformations to get back surfaces applicable to the
initial totally real region.

In what concerns the solitons of quadrics we take as model the
fact that the solitons of the sine-Gordon equation (with the
(vacuum) $0$-soliton $\om=0$) correspond at the level of the
geometric picture to the $0$-soliton being the axis of the
tractrix (thus it is a degenerate surface), the $1$-solitons (B
transforms of the $0$-soliton) being the Dini helicoids (which
include the real pseudo-sphere) and thus one can find the
$n$-solitons, $n\ge 2$ by explicit formulae via the {\it Bianchi
Permutability Theorem} (BPT).

Based on this model the $0$-soliton should be a degenerated
surface (curve or point) and one must be able to explicitly
compute the $1$-solitons (B transforms of the $0$-soliton); after
that the $n$-solitons, $n\ge 2$ will be amenable to explicit
computations via the same BPT.

\section{Totally real forms of the sine-Gordon equation and their
solitons}

Consider the classical hyperbolic sine-Gordon equation
\begin{eqnarray}\label{eq:hsg}
\om_{vv}-\om_{uu}=\cos\om\sin\om
\end{eqnarray} in conjunction with real deformations
$x\subset\mathbb{R}^3$ of the pseudo-sphere (it represents the
Gau\ss\ equation in Chebyshev coordinates $(u+v,u-v)$ which are
further asymptotes) and the classical symmetric {\it B\"{a}cklund}
(B) transformation $\om_1=B_{\si_1}(\om_0),\
\om_0=B_{\si_0}(\om_1),\ \si_0=-\si_1\in\mathbb{R}^*$
\begin{eqnarray}\label{eq:backhsg}
\om_{1v}-\om_{0u}=\frac{\si_1\sin(\om_1+\om_0)+\si_1^{-1}\sin(\om_1-\om_0)}{2},\nonumber\\
\om_{1u}-\om_{0v}=\frac{\si_1\sin(\om_1+\om_0)-\si_1^{-1}\sin(\om_1-\om_0)}{2},\
0\leftrightarrow 1
\end{eqnarray}
of its solutions together with its $1$-solitons $\om_0=0,\
\om_1=\pm
2\tan^{-1}e^{\frac{\si_1-\si_1^{-1}}{2}u+\frac{\si_1+\si_1^{-1}}{2}v+c_1},\
c_1\in\mathbb{R}$ and {\it Bianchi Permutability Theorem} (BPT)
\begin{eqnarray}\label{eq:bpt}
\tan\frac{\om_3-\om_0}{2}=\frac{\si_2+\si_1}{\si_2-\si_1}\tan\frac{\om_2-\om_1}{2}
\ \mathrm{for}\ \om_1=B_{\si_1}(\om_0),\
\om_2=B_{\si_2}(\om_0),\nonumber\\
B_{\si_2}\circ B_{\si_1}(\om_0)=B_{\si_2}(\om_1)=\om_3=
B_{\si_1}(\om_2)=B_{\si_1}\circ B_{\si_2}(\om_0).
\end{eqnarray}
Note that (\ref{eq:bpt}) admits the complex conjugate
$\si_2=\bar\si_1\in\mathbb{C}\setminus\mathbb{R},\
\om_2=\bar\om_1\subset\mathbb{C},\ \om_3,\om_0\subset\mathbb{R}$
version; with certain rationality conditions one obtains at the
level of the geometric picture breathers.

Note also that as it was pointed out by Bianchi when he originally
introduced his BPT in 1890 the BPT does not exclude the case
$\si_2=\si_1$ as being the trivial $\om_3=\om_0$, but allows it as
a limiting case $\si_2\rightarrow \si_1$ and an application of
L'Hospital; for example for $2$ solitons if we let
$c_2=c_2(\si_2),\ c_2(\si_1)=c_1,\ c'_2(\si_1)=\frac{c}{\si_1}$,
then
$\tan\frac{\om_3}{2}=(\frac{\si_1+\si_1^{-1}}{2}u+\frac{\si_1-\si_1^{-1}}{2}v+c)
\sin\om_1$ depends on two constants $c,c_1$ besides the spectral
parameter $\si_1$, as expected (the B transformation should
introduce one constant besides the spectral parameter).

In order to assure that further $n$-th iterates $\mathcal{M}_n$ of
the B transformation ({\it moving M\"{o}bius configurations} in
Bianchi's denomination) give the same result independently of the
chosen path of composition of the B transformation we need to
check only for the third iteration. This is to be expected, since
by discretization the B transformation corresponds to the first
derivative, the BPT ($\mathcal{M}_2$) corresponds to the commuting
of second order derivatives (roughly the Gau\ss-Weingarten
equations or equivalently the flat connection form condition) and
the third M\"{o}bius configuration $\mathcal{M}_3$ corresponds to
the commuting of the third order derivatives (roughly the {\it
Gau\ss-Codazzi-Mainardi-Peterson} (GCMP) equations; as we know
there are no conditions beyond the GCMP equations for a surface).
Moreover the BPT with a leg of the Bianchi quadrilateral
infinitesimal precisely describes the B transformation and thus
the BPT encodes all necessary algebraic information needed to
prove the existence of the B transformation and similarly the
third M\"{o}bius configuration encodes all necessary algebraic
information needed to prove the validity of the BPT.

Since the odd $\mathcal{M}_{2n+1},\ n>0$ M\"{o}bius configuration
does not depend on $\om_0$ we can use $\om_0=0$ to get
$\mathcal{M}_3$
$$(e^{i\om_2+i\om_4}-e^{i\om_1+i\om_7})(\frac{\si_2}{\si_3}-\frac{\si_3}{\si_2})
+(e^{i\om_1+i\om_4}-e^{i\om_2+i\om_7})(\frac{\si_3}{\si_1}-\frac{\si_1}{\si_3})
+(e^{i\om_1+i\om_2}-e^{i\om_4+i\om_7})(\frac{\si_1}{\si_2}-\frac{\si_2}{\si_1})=0.$$
Note again that we can have
$\si_2=\bar\si_1\in\mathbb{C}\setminus\mathbb{R},\
\si_3\in\mathbb{R},\ \om_2=\bar\om_1\subset\mathbb{C},\
\om_4,\om_7\subset\mathbb{R}$.

Note also that unlike the BPT the $\mathcal{M}_3$ configuration is
symmetric in all variables (it has the symmetries of a regular
tetrahedron).
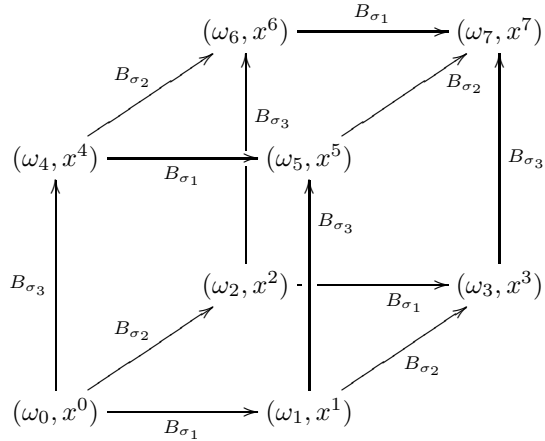
\begin{center}
$\xymatrix@!0{&&&(\om_6,x^6)\ar@{->}[rrrr]^{B_{\si_1}}&&&&(\om_7,x^7)\\
&&&&&&&\\
(\om_4,x^4)\ar@{->}[uurrr]^{B_{\si_2}}
\ar@{->}[rrrr]_{B_{\si_1}}&&&&
(\om_5,x^5)\ar@{->}[uurrr]_>>>>{B_{\si_2}}&&&\\
&&&&&&&\\
&&&(\om_2,x^2)\ar@{->}'[r][rrrr]_{B_{\si_1}}
\ar@{->}'[uu][uuuu]_<<<<<{B_{\si_3}}&&&&
(\om_3,x^3)\ar@{->}[uuuu]_{B_{\si_3}}\\
&&&&&&&\\
(\om_0,x^0)\ar@{->}[rrrr]_{B_{\si_1}}\ar@{->}[uurrr]^{B_{\si_2}}
\ar@{->}[uuuu]^{B_{\si_3}}&&&&(\om_1,x^1)\ar@{->}[uurrr]_{B_{\si_2}}
\ar@{->}[uuuu]_<<<<<<<<<<<<<<<<<<<<<<{B_{\si_3}}&&&}$
\end{center}
Similarly by considering purely imaginary $\om=i\theta$ in the
sine-Gordon equation we have the hyperbolic sinh-Gordon equation
\begin{eqnarray}\label{eq:hsgh}
\theta_{vv}-\theta_{uu}=\cosh\theta\sinh\theta
\end{eqnarray}
with the symmetric B transformation
$\theta_1=B_{\si_1}(\theta_0),\ \theta_0=B_{\si_0}(\theta_1),\
\si_0=-\si_1\in\mathbb{R}^*$
\begin{eqnarray}\label{eq:backhsgh}
\theta_{1v}-\theta_{0u}=\frac{\si_1\sinh(\theta_1+\theta_0)+\si_1^{-1}\sinh(\theta_1-\theta_0)}{2},\nonumber\\
\theta_{1u}-\theta_{0v}=\frac{\si_1\sinh(\theta_1+\theta_0)-\si_1^{-1}\sinh(\theta_1-\theta_0)}{2},\
0\leftrightarrow 1
\end{eqnarray}
of its solution together with its $1$-solitons $\theta_0=0,\
\theta_1=\pm
2\tanh^{-1}e^{\frac{\si_1-\si_1^{-1}}{2}u+\frac{\si_1+\si_1^{-1}}{2}v+c_1},\
c_1\in\mathbb{R},\\ -1<\tanh(\theta_1)<1$, BPT
\begin{eqnarray}\label{eq:bpth}
\tanh\frac{\theta_3-\theta_0}{2}=\frac{\si_2+\si_1}{\si_2-\si_1}\tanh\frac{\theta_2-\theta_1}{2}
\ \mathrm{for}\ \theta_1=B_{\si_1}(\theta_0),\
\theta_2=B_{\si_2}(\theta_0),\nonumber\\
B_{\si_2}\circ B_{\si_1}(\theta_0)=B_{\si_2}(\theta_1)=\theta_3=
B_{\si_1}(\theta_2)=B_{\si_1}\circ B_{\si_2}(\theta_0)
\end{eqnarray}
and third M\"{o}bius configuration $\mathcal{M}_3$
$$(e^{\theta_2+\theta_4}-e^{\theta_1+\theta_7})(\frac{\si_2}{\si_3}-\frac{\si_3}{\si_2})
+(e^{\theta_1+\theta_4}-e^{\theta_2+\theta_7})(\frac{\si_3}{\si_1}-\frac{\si_1}{\si_3})
+(e^{\theta_1+\theta_2}-e^{\theta_4+\theta_7})(\frac{\si_1}{\si_2}-\frac{\si_2}{\si_1})=0.$$
Note that (\ref{eq:bpth}) admits the complex conjugate
$\si_2=\bar\si_1\in\mathbb{C}\setminus\mathbb{R},\
\theta_2=\bar\theta_1\subset\mathbb{C},\
\theta_3,\theta_0\subset\mathbb{R}$ version and we may also have
$\si_2=\bar\si_1\in\mathbb{C}\setminus\mathbb{R},\
\si_3\in\mathbb{R},\ \theta_2=\bar\theta_1\subset\mathbb{C},\
\theta_4,\theta_7\subset\mathbb{R}$.

Similarly by considering purely imaginary coordinate $u$ in the
sine-Gordon equation we have the elliptic sine-Gordon equation
\begin{eqnarray}\label{eq:esg}
\om_{vv}+\om_{uu}=\cos\om\sin\om
\end{eqnarray}
and by considering purely imaginary coordinate $u$  and purely
imaginary $\om=i\theta$ in the sine-Gordon equation we have the
elliptic sinh-Gordon equation
\begin{eqnarray}\label{eq:esgh}
\theta_{vv}+\theta_{uu}=\cosh\theta\sinh\theta
\end{eqnarray}
with the symmetric B transformation $\theta_1=B_{\si_1}(\om_0),\
\om_0=B_{\si_0}(\theta_1),\ \si_0=-\si_1\in\mathbf{S}^1$
\begin{eqnarray}\label{eq:backesgh}
i\theta_{1v}-i\pa_u\om_0=\frac{\si_1\sin(i\theta_1+\om_0)+
\si_1^{-1}\sin(i\theta_1-\om_0)}{2},\nonumber\\
i\pa_u(i\theta_1)-\om_{0v}=\frac{\si_1\sin(i\theta_1+\om_0)
-\si_1^{-1}\sin(i\theta_1-\om_0)}{2},\ 0\leftrightarrow 1
\end{eqnarray}
of their solution together with its $1$-solitons $\om_0=0,\
\theta_1=\pm
2\tanh^{-1}e^{\frac{\si_1-\si_1^{-1}}{2i}u+\frac{\si_1+\si_1^{-1}}{2}v+c_1},\
c_1\in\mathbb{R},\\ -1<\tanh(\theta_1)<1;\ \theta_0=0,\ \om_1=\pm
2\tan^{-1}e^{-\frac{\si_1-\si_1^{-1}}{2i}u-\frac{\si_1+\si_1^{-1}}{2}v+c_1},\
c_1\in\mathbb{R}$, BPT
\begin{eqnarray}\label{eq:bpte}
\tan\frac{\om_3-\om_0}{2}=i\frac{\si_2+\si_1}{\si_2-\si_1}\tanh\frac{\theta_2-\theta_1}{2}
\ \mathrm{for}\ \theta_1=B_{\si_1}(\om_0),\
\theta_2=B_{\si_2}(\om_0),\nonumber\\
B_{\si_2}\circ B_{\si_1}(\om_0)=B_{\si_2}(\theta_1)=\om_3=
B_{\si_1}(\theta_2)=B_{\si_1}\circ B_{\si_2}(\om_0),\nonumber\\
\tanh\frac{\theta_3-\theta_0}{2}=-i\frac{\si_2+\si_1}{\si_2-\si_1}\tan\frac{\om_2-\om_1}{2}
\ \mathrm{for}\ \om_1=B_{\si_1}(\theta_0),\
\om_2=B_{\si_2}(\theta_0),\nonumber\\
B_{\si_2}\circ B_{\si_1}(\theta_0)=B_{\si_2}(\om_1)=\theta_3=
B_{\si_1}(\om_2)=B_{\si_1}\circ B_{\si_2}(\theta_0)
\end{eqnarray}
and third M\"{o}bius configuration $\mathcal{M}_3$
$$(e^{\theta_2+\theta_4}-e^{\theta_1+\theta_7})(\frac{\si_2}{\si_3}-\frac{\si_3}{\si_2})
+(e^{\theta_1+\theta_4}-e^{\theta_2+\theta_7})(\frac{\si_3}{\si_1}-\frac{\si_1}{\si_3})
+(e^{\theta_1+\theta_2}-e^{\theta_4+\theta_7})(\frac{\si_1}{\si_2}-\frac{\si_2}{\si_1})=0,$$
$$(e^{i\om_2+i\om_4}-e^{i\om_1+i\om_7})(\frac{\si_2}{\si_3}-\frac{\si_3}{\si_2})
+(e^{i\om_1+i\om_4}-e^{i\om_2+i\om_7})(\frac{\si_3}{\si_1}-\frac{\si_1}{\si_3})
+(e^{i\om_1+i\om_2}-e^{i\om_4+i\om_7})(\frac{\si_1}{\si_2}-\frac{\si_2}{\si_1})=0.$$
Note that (\ref{eq:bpte}) admits the complex conjugate
$\si_2=\bar\si_1\in\mathbb{C}\setminus\mathbb{R},\
\theta_2=-\bar\theta_1\subset\mathbb{C},\
\om_3,\om_0\subset\mathbb{R};\ \om_2=-\bar\om_1\subset\mathbb{C},\
\theta_3,\theta_0\subset\mathbb{R}$ version and we may also have
$\si_2=\bar\si_1\in\mathbb{C}\setminus\mathbb{R},\
\si_3\in\mathbb{R},\ \theta_2=-\bar\theta_1\subset\mathbb{C},\
\theta_4,\theta_7\subset\mathbb{R};\
\om_2=-\bar\om_1\in\mathbb{C},\ \om_4,\om_7\in\mathbb{R}$.

\section{Bianchi's B\"{a}cklund
transformation for real quadrics via common conjugate systems}

\subsection{Real deformations of (the imaginary region of) the
real hyperbolic paraboloid} \noindent

\noindent

Consider the general confocal real hyperbolic paraboloids in an
isothermic-conjugate parametrization invariant under the Ivory
affinity between confocal quadrics
$$x_z=x_z(\al,\be):=[\sqrt{a_1-z}\al\ \ \sqrt{-a_2+z}\sqrt{\ep}\be\ \ \frac{\al^2-\ep\be^2+z}{2}]^T,\
a_1>z>0>a_2,$$ $$a_1^{-1}-a_2^{-1}=1,\ \al,\be\in\mathbb{R},\
\ep=\pm 1,\ \ep\al^2a_1^{-1}-\be^2a_2^{-1}+\ep>0$$ (the case
$0>z>a_2$ is realized by a rigid motion
$(e_1,e_3)\leftrightarrow(e_2,-e_3)$ and we have imaginary region
for $\ep=-1$) with positive definite linear element, second
fundamental form and Christoffel symbols of $x_0$:
$$|dx_0|^2=a_1d\al^2-\ep a_2d\be^2+(\al d\al-\ep\be d\be)^2;\
N_0^Td^2x_0=\sqrt{\ep}\frac{-\ep d\al^2+d\be^2}{\sqrt{H}},\
H:=\ep\frac{\al^2}{a_1}-\frac{\be^2}{a_2}+\ep;$$
$$-\ep\Ga_{22}^1=\Ga_{11}^1=(\log\sqrt{H})_{\al},\
-\ep\Ga_{11}^2=\Ga_{22}^2=(\log\sqrt{H})_{\be},\
\Ga_{12}^1=\Ga_{12}^2=0.$$ We have the GCMP equations
$$g_{2p}[(\Ga_{11}^p)_2-(\Ga_{12}^p)_1+\Ga_{11}^q\Ga_{q2}^p
-\Ga_{12}^q\Ga_{q1}^p]=R_{1212}=h_{11}h_{22}-h_{12}^2,$$
$$(h_{12})_1-(h_{11})_2+\Ga_{12}^mh_{m1}-\Ga_{11}^mh_{m2}=0,\
(h_{12})_2-(h_{22})_1+\Ga_{12}^mh_{m2}-\Ga_{22}^mh_{m1}=0.$$ We
have a distinguished tangent vector field
$\mathcal{V}_0:=\ep(\log\sqrt{H})_{\al}x_{0\al}-(\log\sqrt{H})_{\be}x_{0\be}$;
it has the properties $|\mathcal{V}_0|^2=1-\frac{1}{H},\
\mathcal{V}_0^Tx_{0\al}=\ep\al,\ \mathcal{V}_0^Tx_{0\be}=-\be$.

Note also that the condition
\begin{eqnarray}\label{eq:pet}
(\Ga_{11}^2\frac{h_{22}}{h_{11}})_{\al}=(\Ga_{22}^1\frac{h_{11}}{h_{22}})_{\be}=
-2\Ga_{11}^2\Ga_{22}^1
\end{eqnarray}
that $(\al,\be)$ is common to an $1$-dimensional Peterson's family
of deformations $x$ of $x_0$ is satisfied.

Given a real deformation $x\subset\mathbb{R}^3$ of a real region
$\subset x_0$ (that is $\ep=1$) there exists a conjugate system
$(u,v)$ common to both $x_0$ and $x$ (this is true for any two
surfaces in a point-wise correspondence). Denote with $\bar\cdot$
the quantities of interest in the GCMP equations (namely the
Christoffel symbols and the second fundamental form) of $x_0$
referred to the $(u,v)$ coordinates and similarly with
$\tilde\cdot$ those of $x$. We have $\al_u\al_v-\be_u\be_v=0$ and
from the Gau\ss\ equation $(\al_u^2-\be_u^2)(\al_v^2-\be_v^2)<0$;
assume (by changing $u$ and $v$ if necessary) $\al_u^2-\be_u^2>0$.
With $\la:=\mathrm{sgn}(\al_u)\sqrt{\al_u^2-\be_u^2},\
\mu:=\mathrm{sgn}(\be_v)\sqrt{\be_v^2-\al_v^2}$ we have $\bar
h_{11}=\frac{\la^2}{\sqrt{H}},\ \bar h_{12}=0,\ \bar
h_{22}=-\frac{\mu^2}{\sqrt{H}}$. From the general formula for the
change of Christoffel symbols $\frac{\pa u^l}{\pa\ti
u^c}\ti\Ga_{ab}^c=\frac{\pa^2u^l}{\pa\ti u^a\pa\ti u^b}+\frac{\pa
u^j}{\pa\ti u^a}\frac{\pa u^k}{\pa\ti u^b}\Ga_{jk}^l$ we get
$\bar\Ga_{12}^1\al_u+\bar\Ga_{12}^2\al_v=\al_{uv},\
\bar\Ga_{12}^1\be_u+\bar\Ga_{12}^2\be_v=\be_{uv},\
\bar\Ga_{11}^1\al_u+\bar\Ga_{11}^2\al_v=\al_{uu}+\la^2(\log\sqrt{H})_{\al},\
\bar\Ga_{11}^1\be_u+\bar\Ga_{11}^2\be_v=\be_{uu}-\la^2(\log\sqrt{H})_{\be},\
\bar\Ga_{22}^1\al_u+\bar\Ga_{22}^2\al_v=\al_{vv}-\mu^2(\log\sqrt{H})_{\al},\
\bar\Ga_{22}^1\be_u+\bar\Ga_{22}^2\be_v=\be_{vv}+\mu^2(\log\sqrt{H})_{\be}$,
so $\bar\Ga_{12}^1=(\log\la)_v,\ \bar\Ga_{12}^2=(\log\mu)_u,\
\bar\Ga_{11}^1=(\log(\la\sqrt{H}))_u,\
\bar\Ga_{11}^2=\frac{\la^2}{\mu^2}(\log\frac{\la}{\sqrt{H}})_v,\
\bar\Ga_{22}^2=(\log(\mu\sqrt{H}))_v,\
\bar\Ga_{22}^1=\frac{\mu^2}{\la^2}(\log\frac{\mu}{\sqrt{H}})_u$.
From the CMP equations of $x_0,\ x$ we have
$$(\bar h_{11})_v=\bar\Ga_{12}^1\bar h_{11}-\bar\Ga_{11}^2\bar
h_{22},\ (\ti h_{11})_v=\bar\Ga_{12}^1\ti h_{11}-\bar\Ga_{11}^2\ti
h_{22},$$
$$(\bar h_{22})_u=\bar\Ga_{12}^2\bar h_{22}-\bar\Ga_{22}^1\bar
h_{11},\ (\ti h_{22})_u=\bar\Ga_{12}^2\ti h_{22}-\bar\Ga_{22}^1\ti
h_{11}.$$ Keeping account of the Gau\ss\ equation $\bar h_{11}\bar
h_{22}=\ti h_{11}\ti h_{22}$ one can multiply the first equations
respectively with $\bar h_{11},\ \ti h_{22}$ (and the second
equations respectively with $\bar h_{22},\ \ti h_{11}$), subtract
them and get rid respectively of the $\bar\Ga_{11}^2,\
\bar\Ga_{22}^1$ terms: $(\log(\bar h_{11}^2-\ti
h_{11}^2))_v=2(\log\la)_v,\ (\log(\bar h_{22}^2-\ti
h_{22}^2))_u=2(\log\mu)_u$. Thus $\bar h_{11}^2-\ti
h_{11}^2=\phi(u)\la^2,\ \bar h_{22}^2-\ti
h_{22}^2=\varphi(v)\mu^2$; after a change of the $u$ and $v$
variables one can absorb $\phi(u),\ \varphi(v)$ up to opposite
signs $\ep_1:=\pm 1,\ \ep_2=-\ep_1$ (here we have again from the
Gau\ss\ equation $\bar h_{11}^2>\ti h_{11}^2\Leftrightarrow\bar
h_{22}^2<\ti h_{22}^2$). We have $\ti h_{11}^2=\bar
h_{11}^2-\ep_1\la^2=\la^2(\frac{\la^2}{H}-\ep_1),\ \ti
h_{22}^2=\bar h_{22}^2+\ep_1\mu^2=\mu^2(\frac{\mu^2}{H}+\ep_1)$;
from the Gau\ss\ equation $\ti h_{11}\ti h_{22}=\bar h_{11}\bar
h_{22}$ we get $H=-\ep_1\mu^2+\ep_1\la^2$; by performing, if
necessary, the change $(\al,u)\leftrightarrow(\be,v)$ we can
choose $\ep_1:=-1$. Thus the second fundamental form of $x$ is
$\frac{\la\mu(du^2-dv^2)}{\sqrt{H}},\ H=\mu^2-\la^2$ and we are
led to consider the hyperbolic angle $\theta$ between the
conjugate systems $(\al,\be)$ and $(u,v)$, that is
$\begin{bmatrix}\al_u&\al_v\\\be_u&\be_v\end{bmatrix}=
\begin{bmatrix}\la\mathbf{C}&\mu\mathbf{S}\\\la\mathbf{S}&\mu\mathbf{C}\end{bmatrix},\
\mathbf{C}:=\cosh\theta,\ \mathbf{S}:=\sinh\theta$ (note that by
doing this the sign of $\theta$ is decided by that of $\be_u$).
Imposing the compatibility conditions
$(\la\mathbf{C})_v=(\mu\mathbf{S})_u,\
(\la\mathbf{S})_v=(\mu\mathbf{C})_u$ we get $\mu_u=\la\theta_v,\
\la_v=\mu\theta_u$; differentiating $H=\mu^2-\la^2$ with respect
to $u$, respectively $v$ we are led to consider the hyperbolic
sinh-Gordon equation (\ref{eq:hsgh}) as the compatibility
condition of the completely integrable linear system in $\al,\
\be,\ \la,\ \mu$:
\begin{eqnarray}\label{eq:linsyst}
d\begin{bmatrix}\al\\\be\\\la\\\mu\end{bmatrix}=\begin{bmatrix}
\la\mathbf{C}du+\mu\mathbf{S}dv\\\la\mathbf{S}du+\mu\mathbf{C}dv\\
(-\mathbf{C}\frac{\al}{a_1}+\mathbf{S}\frac{\be}{a_2}+\mu\theta_v)du+\mu\theta_udv\\
\la\theta_vdu+(\mathbf{S}\frac{\al}{a_1}-\mathbf{C}\frac{\be}{a_2}+\la\theta_u)dv\end{bmatrix},\
\mu^2-\la^2=H.
\end{eqnarray}
Note that a solution $\theta$ of (\ref{eq:hsgh}) will produce an
$1$-dimensional family of deformations $x$ of $x_0$ (from the
original $4$-dimensional space of solutions of the differential
part of (\ref{eq:linsyst}) the prime integral property
$\mu^2-\la^2=H$ removes a constant and translation in $u,v$
another two). The condition that an $1$-dimensional family of
deformations $x$ of $x_0$ with common conjugate system is of
Peterson's type (that is (\ref{eq:pet}) is satisfied in the
$(u,v)$ coordinates) is invariant under changes of variables
$(u,v)$ into themselves; in our case we need
$(\log\frac{\la}{\mu})_{uv}=0$, but this adjoined to
(\ref{eq:linsyst}) will be over-determined; as we shall see later
the condition $(\log\frac{\la}{\mu})_{uv}=0$ will be preserved by
the B transformation.

Note that if we assume that the common conjugate system on $x_0,\
x$ is isothermic-conjugate on $x$ (Darboux), then from the Gau\ss\
equations we obtain immediately that the second fundamental form
of $x$ is $\frac{\la\mu(du^2-dv^2)}{\sqrt{H}}$; everything else
except $\mu_v,\ \la_u,\ H=\mu^2-\la^2$ follows immediately as
previously. The remaining needed information follows immediately
from the CMP equations of $x$:
$d(\log\frac{\la\mu}{\sqrt{H}})=(\bar\Ga_{12}^2+\bar\Ga_{22}^1)du+(\bar\Ga_{12}^1+\bar\Ga_{11}^2)dv$
become $d\log\frac{\mu^2-\la^2}{H}=0$; by a same homothety in the
$(u,v)$ variables and a choice of sign we can assume
$\frac{\mu^2-\la^2}{H}=1$.

By similar computations if $x\subset\mathbb{R}^3$ is a real
deformation of an imaginary region $\subset x_0$ (that is
$\ep=-1$), then with $\la:=\pm\sqrt{\al_u^2+\be_u^2},\
\mu:=\pm\sqrt{\al_v^2+\be_v^2}$ (the signs may vary when we shall
consider the B transformation) the second fundamental form of $x$
is $\frac{\la\mu(du^2-dv^2)}{\sqrt{H}},\ H=\mu^2+\la^2$ and we are
led to consider the angle $\om$ between the conjugate systems
$(\al,\be)$ and $(u,v)$, that is
$\begin{bmatrix}\al_u&\al_v\\\be_u&\be_v\end{bmatrix}=
\begin{bmatrix}\la\mathbf{C}&-\mu\mathbf{S}\\\la\mathbf{S}&\mu\mathbf{C}\end{bmatrix},\
\mathbf{C}:=\cos\om,\ \mathbf{S}:=\sin\om$. Imposing the
compatibility conditions $(\la\mathbf{C})_v=-(\mu\mathbf{S})_u,\
(\la\mathbf{S})_v=(\mu\mathbf{C})_u$ we get $\mu_u=\la\om_v,\
\la_v=-\mu\om_u$; differentiating $H=\mu^2+\la^2$ with respect to
$u$, respectively $v$ we are led to consider the hyperbolic
sine-Gordon equation (\ref{eq:hsg}) as the compatibility condition
of the completely integrable linear system in $\al,\ \be,\ \la,\
\mu$:
\begin{eqnarray}\label{eq:linsysti}
d\begin{bmatrix}\al\\\be\\\la\\\mu\end{bmatrix}=\begin{bmatrix}
\la\mathbf{C}du-\mu\mathbf{S}dv\\\la\mathbf{S}du+\mu\mathbf{C}dv\\
(-\mathbf{C}\frac{\al}{a_1}-\mathbf{S}\frac{\be}{a_2}-\mu\om_v)du-\mu\om_udv\\
\la\om_vdu+(\mathbf{S}\frac{\al}{a_1}-\mathbf{C}\frac{\be}{a_2}+\la\om_u)dv\end{bmatrix},\
\mu^2+\la^2=H.
\end{eqnarray}
One can put everything in  a matrix notation: with
$R:=\begin{bmatrix}\mathbf{C}&\ep\mathbf{S}\\\mathbf{S}&\mathbf{C}\end{bmatrix},\
\del:=\mathrm{diag}[du\ \ dv],\ V:=[\al\ \ \be]^T,\ \La:=[\la\ \
\mu]^T,\ A':=\mathrm{diag}[a_1^{-1}\ \ a_2^{-1}],\
\Om:=R^{-1}R_udv+R^{-1}R_vdu,\ \mathcal{E}:=\mathrm{diag}[1\ \
-\ep]$ we have the sine(sinh)-Gordon equation
$$d\wedge\Om-\Om\wedge\Om=-\del R^{-1}A'\wedge R\del,\
R^{-1}dR\wedge \del-\del\wedge\Om=0 \Leftrightarrow$$
$$e_1^T[(R^{-1}R_u)_u-(R^{-1}R_v)_v+R^{-1}A'R]e_2=0$$ as the compatibility condition for the
completely integrable linear system
$$d\begin{bmatrix}V\\\La\end{bmatrix}=\begin{bmatrix}0&R\del\\
-\del
R^{-1}A'&\Om\end{bmatrix}\begin{bmatrix}V\\\La\end{bmatrix},\
\La^T\mathcal{E}\La=-V^T\mathcal{E}A'V-1.$$ With $R'_z:=I_2-zA'$
consider two points $x_0^0,\ x_0^1\in x_0$ in the same totally
real region of $x_0$ (corresponding to $\ep=\pm 1$) such that
$x_0^0,\ x_z^1$ are in the symmetric tangency configuration
\begin{eqnarray}\label{eq:tc}
x_z^1\in T_{x_0^0}x_0\Leftrightarrow x_z^1=x_0^0+[x_{0\al_0}^0\ \
x_{0\be_0}^0](\sqrt{R'_z}V_1-V_0)\Leftrightarrow\nonumber\\
(\sqrt{R'_z}V_1-V_0)^T\mathcal{E}(\sqrt{R'_z}V_1-V_0)=-\ep zH_1,\
0\leftrightarrow 1.
\end{eqnarray}
Thus a quadratic functional relationship is established between
$\al_0,\ \be_0,\ \al_1,\ \be_1$ and only three among them remain
functionally independent:
\begin{eqnarray}\label{eq:dtc}
dV_0^T\mathcal{E}(\sqrt{R'_z}V_1-V_0)=-dV_1^T\mathcal{E}(\sqrt{R'_z}V_0-V_1).
\end{eqnarray}
Bianchi's main theorem on the theory of deformations of quadrics
states (in our case) that given the {\it seed} deformation
$x^0\subset\mathbb{R}^3$ of the totally real region $x_0^0\subset
x_0$ (that is $|dx^0|^2=|dx_0^0|^2$) the differential system
\begin{eqnarray}\label{eq:weingp}
x^1=x^0+[x_{\al_0}^0\ \ x_{\be_0}^0](\sqrt{R'_z}V_1-V_0),\
|dx^1|^2=|dx_0^1|^2
\end{eqnarray}
obtained by imposing the ACPIA $|dx^1|^2=|dx_0^1|^2$ is completely
integrable (thus it admits two $1$-dimensional family of solutions
({\it leaves}) $x^1=:B_{z,\ep_1}(x^0)\subset\mathbb{R}^3,\
\ep_1=\pm 1$ whose determination requires the integration of a
Ricatti equation), that $x^0,x^1$ are the focal surfaces of a {\it
Weingarten congruence} (congruence of lines on whose two focal
surfaces the asymptotic directions correspond; since conjugate
directions are harmonically conjugate to the asymptotic ones all
conjugate systems correspond in this case), that $x^1$ is
applicable to $x_0^1$ (in our case of the same totally real region
of $x_0$ as $x^0_0$) and that we have the symmetry
$0\leftrightarrow 1$. Its simplest proof uses parametrization by
rulings, since they behave well with respect to the metric
properties of the Ivory affinity and for the particular
configuration $x^0=x_0^0$ the leaves $x^1$ become rulings on
$x_z$; it appears elsewhere so we shall not insist on it here.

Using (\ref{eq:dtc}) we get by differentiating (\ref{eq:weingp})
$$dx^1=-\ep dV_1^T\mathcal{E}(\sqrt{R'_z}V_0-V_1)\mathcal{V}^0+[x_{\al_0}^0\ \
x_{\be_0}^0]\sqrt{R'_z}dV_1$$$$+N^0(N^0)^T[dx_{\al_0}^0\ \
dx_{\be_0}^0](\sqrt{R'_z}V_1-V_0),\ N^0:=\frac{x_{\al_0}^0\times
x_{\be_0}^0}{\sqrt{-a_1a_2H_0}},$$ so the ACPIA becomes
\begin{eqnarray}\label{eq:dx1}
(N^0)^T[dx_{\al_0}^0\ \
dx_{\be_0}^0](\sqrt{R'_z}V_1-V_0)=\frac{\ep_1}{\sqrt{H_0}}dV_1^T
\begin{bmatrix}0&-1\\1&0\end{bmatrix}(\sqrt{R'_z}V_0-V_1),\nonumber\\
\ep_0=\ep_1=\pm 1,\ 0\leftrightarrow 1.
\end{eqnarray}
Using
$\begin{bmatrix}\frac{\pa}{\pa\al_0}\\\frac{\pa}{\pa\be_0}\end{bmatrix}=
\begin{bmatrix}\frac{\mathbf{C}_0}{\la_0}&-\frac{\mathbf{S}_0}{\mu_0}\\
-\ep\frac{\mathbf{S}_0}{\la_0}&\frac{\mathbf{C}_0}{\mu_0}\end{bmatrix}
\begin{bmatrix}\frac{\pa}{\pa
u}\\\frac{\pa}{\pa v}\end{bmatrix}$, (\ref{eq:dtc}) and
(\ref{eq:dx1}) we get with $r_j:=\sqrt{1-za_j^{-1}},\ j=1,2$:
$\\\mu_0[(r_1\al_1-\al_0)\mathbf{C_0}-\ep(r_2\be_1-\be_0)\mathbf{S_0}]=
\ep_1[(r_2\be_0-\be_1)\al_{1u}-(r_1\al_0-\al_1)\be_{1u}],\\
\la_0[(r_2\be_1-\be_0)\mathbf{C_0}-(r_1\al_1-\al_0)\mathbf{S_0}]=
\ep_1[(r_1\al_0-\al_1)\be_{1v}-(r_2\be_0-\be_1)\al_{1v}],\\
\la_0[(r_1\al_1-\al_0)\mathbf{C_0}-\ep(r_2\be_1-\be_0)\mathbf{S_0}]=
-[(r_1\al_0-\al_1)\al_{1u}-\ep(r_2\be_0-\be_1)\be_{1u}],\\
\mu_0[(r_2\be_1-\be_0)\mathbf{C_0}-(r_1\al_1-\al_0)\mathbf{S_0}]=
-[(r_2\be_0-\be_1)\be_{1v}-\ep(r_1\al_0-\al_1)\al_{1v}]$.

If these equations are $I-IV$ and using
$(r_2\be_1-\be_0)^2-\ep(r_1\al_1-\al_0)^2=zH_1>0$, then by
considering $I^2-\ep III^2,\ IV^2-\ep II^2,\ I\cdot II-III\cdot
IV$ we obtain $\al_{1u}^2-\ep\be_{1u}^2>0,\
\be_{1v}^2-\ep\al_{1v}^2>0,\\
\al_{1u}\al_{1v}-\ep\be_{1u}\be_{1v}=0$, so
$\begin{bmatrix}\al_{1u}&\al_{1v}\\\be_{1u}&\be_{1v}\end{bmatrix}=
\begin{bmatrix}\la_1\mathbf{C}_1&\ep\mu_1\mathbf{S}_1\\
\la_1\mathbf{S}_1&\mu_1\mathbf{C}_1\end{bmatrix}$. Because of the
symmetry $0\leftrightarrow 1$ and using $I,\ II$ we obtain that
the second fundamental form of $x^1$ is
$\frac{\la_1\mu_1(du^2-dv^2)}{\sqrt{H_1}}$, so $(u,v)$ is also
isothermic-conjugate on $x^1$ and the B transformation preserves
the orientation of $(u,v)$; from $II^2-\ep III^2$ we obtain
$\mu_1^2-\ep\la_1^2=H_1$ and we have complete symmetry
$0\leftrightarrow 1$ also at the level of the isothermic-conjugate
system $(u,v)$. Thus
$\\\al_1=r_1\al_0-\ep_1'\sqrt{z}(\mathbf{C}_1\la_0-\ep\ep_1\mathbf{S}_1\mu_0),\
\be_1=r_2\be_0-\ep_1'\sqrt{z}(\mathbf{S}_1\la_0-\ep_1\mathbf{C}_1\mu_0),\
\ep_0'=-\ep_1',\ 0\leftrightarrow 1$, or
\begin{eqnarray}\label{eq:algsoli}
\begin{bmatrix}V_1\\\La_1\end{bmatrix}=\ep'_1\sqrt{z}\begin{bmatrix}
I_2&0\\0&\mathcal{E}_1R_0^{-1}\end{bmatrix}\begin{bmatrix}
D&-I_2\\
A'&D\end{bmatrix}\begin{bmatrix}
I_2&0\\0&R_1\mathcal{E}_1\end{bmatrix}\begin{bmatrix}V_0\\\La_0\end{bmatrix},\
D:=\frac{\sqrt{R'_z}}{\ep'_1\sqrt{z}},\ 0\leftrightarrow 1.
\end{eqnarray}
Finally differentiating $V_1$ we obtain the B transformation at
the analytic level
\begin{eqnarray}\label{eq:backl}
dR_1\mathcal{E}_1=-R_1\mathcal{E}_1\Om_0-R_1\mathcal{E}_1\del
R_0^{-1}DR_1\mathcal{E}_1+DR_0\del;
\end{eqnarray}
with $\mathrm{diag}[\frac{\si_1-\si_1^{-1}}{2}\ \
\frac{\si_1+\si_1^{-1}}{2}]:=-D\mathcal{E}_1$ this becomes the B
transformations (\ref{eq:backhsg}) respectively
(\ref{eq:backhsgh}) between solutions $\om_0,\ \ep_1\om_1$ of
(\ref{eq:hsg}) and respectively $\theta_0,\ \ep_1\theta_1$ of
(\ref{eq:hsgh}).

The BPT states that if $x^j=B_{z_j,\ep_j}(x^0),\ j=1,2$, then one
can find only by algebraic computations and derivatives a surface
$x^3$ such that $B_{z_2,\ep_2}(x^1)=x^3=B_{z_1,\ep_1}(x^2)$, that
is $B_{z_1,\ep_1}\circ B_{z_2,\ep_2}=B_{z_2,\ep_2}\circ
B_{z_1,\ep_1}$. Again one can derive the analytic BPT for the B
transformation of the hyperbolic sine(sinh)-Gordon equation from
the geometric picture, just as we did for the B transformation
itself, but we take advantage of the already completed algebraic
transformation of solutions to derive the BPT at the analytic
level, following that we shall then use these analytic
computations to get the geometric realization of solutions in
space.

We have $\begin{bmatrix}
I_2&0\\0&\mathcal{E}_2R_1^{-1}\end{bmatrix}\begin{bmatrix}
D_2&-I_2\\
A'&D_2\end{bmatrix}\begin{bmatrix}
I_2&0\\0&R_3\mathcal{E}_2\mathcal{E}_1R_0^{-1}\end{bmatrix}\begin{bmatrix}
D_1&-I_2\\
A'&D_1\end{bmatrix}\begin{bmatrix}
I_2&0\\0&R_1\mathcal{E}_1\end{bmatrix}\\
=\begin{bmatrix}
I_2&0\\0&\mathcal{E}_1R_2^{-1}\end{bmatrix}\begin{bmatrix}
D_1&-I_2\\
A'&D_1\end{bmatrix}\begin{bmatrix}
I_2&0\\0&R_3\mathcal{E}_1\mathcal{E}_2R_0^{-1}\end{bmatrix}\begin{bmatrix}
D_2&-I_2\\
A'&D_2\end{bmatrix}\begin{bmatrix}
I_2&0\\0&R_2\mathcal{E}_2\end{bmatrix}$, or
$$R_3\mathcal{E}_1\mathcal{E}_2R_0^{-1}=(D_1R_2\mathcal{E}_1\mathcal{E}_2R_1^{-1}-D_2)
(D_1-D_2R_2\mathcal{E}_1\mathcal{E}_2R_1^{-1})^{-1};$$ again this
is just (\ref{eq:bpt}) and (\ref{eq:bpth}) and for
$\ep_1\ep_2\om_3$ (resp $\ep_1\ep_2\theta_3$) we have only two
cases $\ep_1=\pm\ep_2$ in what concerns the dependence on
$\ep_1,\ep_2$ because the BPT requires a $\mathbb{Z}_2$ co-cycle
condition on these signs. We have the space realization

$x^3=x^0-\ep(\sqrt{R'_{z_2}}V_3-V_1)^T\mathcal{E}(\sqrt{R'_{z_1}}V_0-V_1)\mathcal{V}_0+
[x_{\al_0}^0\ \
x_{\be_0}^0](\sqrt{R'_{z_1}}\sqrt{R'_{z_2}}V_3-V_0)+\frac{\ep_1}{\sqrt{H_0}}
(\sqrt{R'_{z_2}}V_3-V_1)^T\begin{bmatrix}0&-1\\1&0\end{bmatrix}(\sqrt{R'_{z_1}}V_0-V_1)N^0,\
(V_1,z_1,\ep_1)\leftrightarrow(V_2,z_2,\ep_2)$.

For the $\mathcal{M}_3$ configuration consider
$(D_1D_2)^{-1}[(D_1^2-D_2^2)
D_1R_1\mathcal{E}_1(D_1R_1\mathcal{E}_1-D_2R_2\mathcal{E}_2)^{-1}-D_1^2]
=R_3\mathcal{E}_1\mathcal{E}_2R_0^{-1}
=(D_1D_2)^{-1}[(D_1^2-D_2^2)
D_2R_2\mathcal{E}_2(D_1R_1\mathcal{E}_1-D_2R_2\mathcal{E}_2)^{-1}-D_2^2],\
R_5\mathcal{E}_1\mathcal{E}_3R_0^{-1}=(D_1D_3)^{-1}[\\(D_3^2-D_1^2)
D_1R_1\mathcal{E}_1(D_3R_4\mathcal{E}_3-D_1R_1\mathcal{E}_1)^{-1}-D_1^2],\
R_6\mathcal{E}_2\mathcal{E}_3R_0^{-1}=(D_2D_3)^{-1}[(D_2^2-D_3^2)
D_2R_2\mathcal{E}_2(D_2R_2\mathcal{E}_2\\-D_3R_4\mathcal{E}_3)^{-1}-D_2^2]$;
thus with
$\Box:=(D_2^2-D_3^2)D_1R_1\mathcal{E}_1+(D_3^2-D_1^2)D_2R_2\mathcal{E}_2+
(D_1^2-D_2^2)D_3R_4\mathcal{E}_3$ we have
$(D_2R_3\mathcal{E}_1\mathcal{E}_2R_0^{-1}-D_3R_5\mathcal{E}_1\mathcal{E}_3R_0^{-1})^{-1}R_1\mathcal{E}_1
=[(D_1^2-D_2^2)(D_1R_1\mathcal{E}_1-D_2R_2\mathcal{E}_2)^{-1}-
(D_3^2-D_1^2)(D_3R_4\mathcal{E}_3-D_1R_1\mathcal{E}_1)^{-1}]^{-1}
=(D_1R_1\mathcal{E}_1-D_2R_2\mathcal{E}_2)\Box^{-1}(D_3R_4\mathcal{E}_3
-D_1R_1\mathcal{E}_1)$ and similarly
$\\(D_3R_6\mathcal{E}_2\mathcal{E}_3R_0^{-1}-D_1R_3\mathcal{E}_1\mathcal{E}_2R_0^{-1})^{-1}R_2\mathcal{E}_2=
(D_1R_1\mathcal{E}_1-D_2R_2\mathcal{E}_2)\Box^{-1}
(D_2R_2\mathcal{E}_2-D_3R_4\mathcal{E}_3)$. Now
$D_1[(D_2^2-D_3^2)D_2R_3\mathcal{E}_1\mathcal{E}_2R_0^{-1}
(D_2R_3\mathcal{E}_1\mathcal{E}_2R_0^{-1}-D_3R_5\mathcal{E}_1\mathcal{E}_3R_0^{-1})^{-1}
R_1\mathcal{E}_1-D_2^2R_1\mathcal{E}_1]=
(D_1^2D_2R_2\mathcal{E}_2-D_2^2D_1R_1\mathcal{E}_1)\Box^{-1}(D_2^2-D_3^2)
(D_3R_4\mathcal{E}_3-D_1R_1\mathcal{E}_1)-D_2^2D_1R_1\mathcal{E}_1
=(D_1^2D_2R_2\mathcal{E}_2-D_2^2D_1R_1\mathcal{E}_1)\Box^{-1}(D_3^2-D_1^2)
(D_2R_2\mathcal{E}_2-D_3R_4\mathcal{E}_3)
-D_1^2D_2R_2\mathcal{E}_2
=D_2[(D_3^2-D_1^2)D_1R_3\mathcal{E}_1\mathcal{E}_2R_0^{-1}
(D_3R_6\mathcal{E}_2\mathcal{E}_3R_0^{-1}-\\
D_1R_3\mathcal{E}_1\mathcal{E}_2R_0^{-1})^{-1}R_2\mathcal{E}_2
-D_1^2R_2\mathcal{E}_2]$, so the very {\it left hand side} (lhs)
and {\it right hand side} (rhs) provide the good definition of and
afford themselves the name
$D_1D_2D_3R_7\mathcal{E}_1\mathcal{E}_2\mathcal{E}_3$. Again these
correspond to the $\mathcal{M}_3$ configurations of the
sine(sinh)-Gordon equations.

Finally to get the first three iterations of the B transformation
for Peterson's real deformations of totally real regions of the
real hyperbolic paraboloid it is enough to give only
$\al_0,\be_0,\la_0,\mu_0$ and the space realization of $x^0$;
everything else will follow according to the established algebraic
formulae.

For $\ep=1$ and $\theta_0=0$ we have $\al_0=\al_0(u),\
\la_0=\la_0(u)=\al_0'(u),\ \be_0=\be_0(v),\
\mu_0=\mu_0(v)=\be_0'(v),\
\mu_0^2+\frac{\be_0^2}{a_2}=\la_0^2+\frac{\al_0^2}{a_1}+1$, so
$\al_0=\sqrt{a_1}\sinh s\sin\frac{u}{\sqrt{a_1}},\
\be_0=\sqrt{-a_2}\cosh s\sinh\frac{v}{\sqrt{-a_2}}$ and we have
Peterson's $1$-dimensional family
$x^0(\al_0,\be_0,s):=[\int_0^{\al_0}\sqrt{a_1-\frac{t^2}{\sinh^2s}}dt\
\ \int_0^{\be_0}\sqrt{-a_2+\frac{t^2}{\cosh^2s}}dt\\
\frac{\al_0^2-\be_0^2\tanh^2s}{2\tanh
s}]^T=\frac{1}{4}[a_1(\frac{2u}{\sqrt{a_1}}+\sin\frac{2u}{\sqrt{a_1}})\sinh
s\ \ -a_2(\frac{2v}{\sqrt{-a_2}}+\sinh\frac{2v}{\sqrt{-a_2}})\cosh
s\ \
(a_1\sin^2\frac{u}{\sqrt{a_1}}+\\a_2\sinh^2\frac{v}{\sqrt{-a_2}})\sinh
2s]^T,\ s\in(0,\infty)$ of real deformations of
$x_0^0=x^0(\al_0,\be_0,\infty)$ which preserves a conjugate
system.

For $\ep=-1$ and $\om_0=0$ we have $\al_0=\al_0(u),\
\la_0=\la_0(u)=\al_0'(u),\ \be_0=\be_0(v),\
\mu_0=\mu_0(v)=\be_0'(v),\
-\mu_0^2-\frac{\be_0^2}{a_2}=\la_0^2+\frac{\al_0^2}{a_1}+1$, so
$\al_0=\sqrt{a_1}\sinh s\sin\frac{u}{\sqrt{a_1}},\
\be_0=\sqrt{-a_2}\cosh
s\cosh\frac{v}{\sqrt{-a_2}}>\sqrt{-a_2}\cosh s,\ v>0$ and we have
Peterson's $1$-dimensional family
$x^0(\al_0,\be_0,s):=[\int_0^{\al_0}\sqrt{a_1-\frac{t^2}{\sinh^2s}}dt\\
\int_{\sqrt{-a_2}\cosh
s}^{\be_0}\sqrt{a_2+\frac{t^2}{\cosh^2s}}dt\ \
\frac{\al_0^2+\be_0^2\tanh^2s}{2\tanh
s}]^T=\frac{1}{4}[a_1(\frac{2u}{\sqrt{a_1}}+\sin\frac{2u}{\sqrt{a_1}})\sinh
s\ \ a_2(\frac{2v}{\sqrt{-a_2}}-\sinh\frac{2v}{\sqrt{-a_2}})\cosh
s\\
(a_1\sin^2\frac{u}{\sqrt{a_1}}-a_2\sinh^2\frac{v}{\sqrt{-a_2}})\sinh
2s]^T,\ s\in(0,\infty)$ of real deformations of
$x_0^0=[\sqrt{a_1}\al_0\ \ \sqrt{-a_2}i\be_0\\
\frac{\al_0^2+\be_0^2}{2}]^T$ which preserves a conjugate system.

Note that one can begin with different Peterson's deformations as
seed: we take the pencil of planes passing through the $e_3$ axis;
the other pencil of planes will be planes perpendicular on the
$e_3$ axis.

Thus for $\ep=1$ we have
$x_0=x_0(\ti\al,\ti\be):=[\sqrt{a_1}e^{\ti\al}\cosh\ti\be\ \
\sqrt{-a_2}e^{\ti\al}\sinh\ti\be\ \ \frac{e^{2\ti\al}}{2}]^T$ with
linear element, second fundamental form and Christoffel symbols
$\\|dx_0|^2=e^{2\ti\al}[(a_1-a_2)(\cosh\ti\be d\ti\al+\sinh\ti\be
d\ti\be)^2+a_2(d\ti\al^2-d\ti\be^2)+e^{2\ti\al}d\ti\al^2],\
N_0^Td^2x_0=\frac{-d\ti\al^2+d\ti\be^2}{e^{-\ti\al}\sqrt{\ti H}},\
\ti H:=\sinh^2\ti\be+a_1^{-1}+e^{-2\ti\al},\ \ti\Ga_{12}^1=0,\
\ti\Ga_{12}^2=1,\ 2-\ti\Ga_{11}^1=\ti\Ga_{22}^1=-(\log\sqrt{\ti
H})_{\ti\al},\ -\ti\Ga_{11}^2=\ti\Ga_{22}^2=(\log\sqrt{\ti
H})_{\ti\be}$.

Again with
$\ti\la:=\mathrm{sgn}(\ti\al_u)\sqrt{\ti\al_u^2-\ti\be_u^2},\
\ti\mu:=\mathrm{sgn}(\ti\be_v)\sqrt{\ti\be_v^2-\ti\al_v^2},\
\begin{bmatrix}\ti\al_u&\ti\al_v\\\ti\be_u&\ti\be_v\end{bmatrix}=
\begin{bmatrix}\ti\la\ti{\mathbf{C}}&\ti\mu\ti{\mathbf{S}}\\
\ti\la\ti{\mathbf{S}}&\ti\mu\ti{\mathbf{C}}\end{bmatrix},\
\ti{\mathbf{C}}:=\cosh\ti\theta,\ \ti{\mathbf{S}}:=\sinh\ti\theta$
we have $(\mu^2,\la^2,H)=e^{2\ti\al}(\ti\mu^2,\ti\la^2,\ti H)$, so
$\ti\mu^2-\ti\la^2=\ti H,\ \bar
h_{11}=\frac{\ti\la^2}{e^{-\ti\al}\sqrt{\ti H}},\ \bar h_{12}=0,\
\bar h_{22}=-\frac{\ti\mu^2}{e^{-\ti\al}\sqrt{\ti H}},\
\bar\Ga_{12}^1=(\log e^{\ti\al}\ti\la)_v,\ \bar\Ga_{12}^2=(\log
e^{\ti\al}\ti\mu)_u,\
\bar\Ga_{11}^2=\frac{\ti\la^2}{\ti\mu^2}(\log\frac{\ti\la}{\sqrt{\ti
H}})_v,\
\bar\Ga_{22}^1=\frac{\ti\mu^2}{\ti\la^2}(\log\frac{\ti\mu}{\sqrt{\ti
H}})_u$.

Thus we get the differential system
\begin{eqnarray}\label{eq:diffsyst}
d\begin{bmatrix}\ti\al\\\ti\be\\\ti\la\\\ti\mu\end{bmatrix}=\begin{bmatrix}
\ti\la\ti{\mathbf{C}}du+\ti\mu\ti{\mathbf{S}}dv\\\ti\la\ti{\mathbf{S}}du+\ti\mu\ti{\mathbf{C}}dv\\
(\ti{\mathbf{C}}e^{-2\ti\al}-\ti{\mathbf{S}}\sinh\ti\be\cosh\ti\be
+\ti\mu\ti\theta_v)du+\ti\mu\ti\theta_udv\\
\ti\la\ti\theta_vdu+(-\ti{\mathbf{S}}e^{-2\ti\al}+\ti{\mathbf{C}}\sinh\ti\be\cosh\ti\be
+\ti\la\ti\theta_u)dv\end{bmatrix},\ \ti\mu^2-\ti\la^2=\ti H
\end{eqnarray}
with compatibility condition
\begin{eqnarray}\label{eq:mhypsg}
\ti\theta_{vv}-\ti\theta_{uu}=(\cosh
2\ti\be+2e^{-2\ti\al})\ti{\mathbf{S}}\ti{\mathbf{C}}.
\end{eqnarray}
Note that we need further manipulation, since the dependence of
$\ti\al,\ti\be$ on $(u,v)$ in (\ref{eq:mhypsg}) is a-priori
undetermined, but that it is not our interest right now, since we
are interested only in the solution $\ti\theta=0$. In this case we
have $\ti\al=\ti\al(u),\ \ti\la=\ti\la(u)=\ti\al'(u),\
\ti\be=\ti\be(v),\ \ti\mu=\ti\mu(v)=\ti\be'(v),\
\ti\be'^2(v)-\sinh^2\ti\be(v)=\ti\al'^2(u)+e^{-2\ti\al(u)}+a_1^{-1}=c$.

\subsection{Imaginary
deformations of (the imaginary region of) the real hyperbolic
paraboloid} \noindent

\noindent

In this case $x^1$ will be applicable to $x_0^1$ of a different
totally real region of $x_0$ as $x^0_0$, so the confocal $x_z^1$
will change type from a real metric point of view. Thus we are led
to consider the elliptic paraboloids
$$x_z=x_z(\al,\be):=[\sqrt{a_1-z}\al\ \ \sqrt{a_2-z}\sqrt{\ep}\be\ \ \frac{\al^2+\ep\be^2+z}{2}]^T,\
a_2>z$$ (the case $z>a_1$ is realized by a rigid motion
$(e_1,e_3)\leftrightarrow(e_2,-e_3)$ and we have imaginary region
for $\ep=-1$) confocal to the given hyperbolic one $x_0$. The
Ivory affinity $\\x_0\mapsto\mathrm{diag}[\sqrt{1-za_1^{-1}}\ \
\sqrt{1-za_2^{-1}}\ \ 1]x_0+\frac{z}{2}=[\sqrt{a_1-z}\al\ \
\sqrt{a_2-z}i\sqrt{\ep}\be\ \ \frac{\al^2-\ep\be^2+z}{2}]^T$ takes
in this case the real (imaginary) region of $x_0$ to the imaginary
(real) one of $x_z$.

If $x\subset\mathbb{R}^2\times i\mathbb{R}$ is an imaginary
deformation of a totally real region $\subset x_0$, then by
similar computations its second fundamental form is
$i\frac{\la\mu(du^2+dv^2)}{\sqrt{H}},\ H=\mu^2+\ep\la^2$ and we
are led to consider the elliptic sinh-Gordon equation
(\ref{eq:esgh}) as the compatibility condition of the completely
integrable linear system in $\al,\ \be,\ \la,\ \mu$:
\begin{eqnarray}\label{eq:linsyste}
d\begin{bmatrix}\al\\\be\\\la\\\mu\end{bmatrix}=\begin{bmatrix}
\la\mathbf{C}du+\mu\mathbf{S}dv\\\la\mathbf{S}du+\mu\mathbf{C}dv\\
(\mathbf{C}\frac{\al}{a_1}-\mathbf{S}\frac{\be}{a_2}-\mu\theta_v)du+\mu\theta_udv\\
\la\theta_vdu+(\mathbf{S}\frac{\al}{a_1}-\mathbf{C}\frac{\be}{a_2}-\la\theta_u)dv\end{bmatrix},\
\mu^2+\la^2=H
\end{eqnarray}
for $\ep=1$ and the elliptic sine-Gordon equation (\ref{eq:esg})
as the compatibility condition of the completely integrable linear
system in $\al,\ \be,\ \la,\ \mu$:
\begin{eqnarray}\label{eq:linsystie}
d\begin{bmatrix}\al\\\be\\\la\\\mu\end{bmatrix}=\begin{bmatrix}
\la\mathbf{C}du-\mu\mathbf{S}dv\\\la\mathbf{S}du+\mu\mathbf{C}dv\\
(\mathbf{C}\frac{\al}{a_1}+\mathbf{S}\frac{\be}{a_2}+\mu\om_v)du-\mu\om_udv\\
\la\om_vdu+(\mathbf{S}\frac{\al}{a_1}-\mathbf{C}\frac{\be}{a_2}-\la\om_u)dv\end{bmatrix},\
\mu^2-\la^2=H
\end{eqnarray}
for $\ep=-1$.

With $r_1:=\sqrt{1-za_1^{-1}},\ r_2:=\sqrt{-1+za_2^{-1}}$ we have
\begin{eqnarray}\label{eq:weingpi}
x_z^1=x_0^0+(r_1\al_1-\al_0)x_{0\al_0}^0 +(-\ep
r_2\be_1-\be_0)x_{0\be_0}^0,\ x^1=x^0+(r_1\al_1-\al_0)x_{\al_0}^0
+(-\ep r_2\be_1-\be_0)x_{\be_0}^0,\nonumber\\
x^1_{\al_1}=-\ep(r_1\al_0-\al_1)(\frac{\al_0}{a_1H_0}x_{\al_0}^0+
\frac{\be_0}{a_2H_0}x_{\be_0}^0)+r_1x_{\al_0}^0
+i\frac{\ep_1(\ep r_2\be_0-\be_1)}{\sqrt{H_0}}N^0,\nonumber\\
x^1_{\be_1}=-(\ep r_2\be_0-\be_1)(\frac{\al_0}{a_1H_0}x_{\al_0}^0+
\frac{\be_0}{a_2H_0}x_{\be_0}^0)-\ep r_2x_{\be_0}^0
-i\frac{\ep_1(r_1\al_0-\al_1)}{\sqrt{H_0}}N^0,\nonumber\\
N^0:=\frac{x_{\al_0}^0\times
x_{\be_0}^0}{\sqrt{-a_1a_2H_0}}\subset(i\mathbb{R})^2\times\mathbb{R},\
\ep_0=\ep_1=\pm 1,\ (0,\ep)\leftrightarrow(1,-\ep).
\end{eqnarray}
We get $(r_1\al_1-\al_0)(N^0)^Tdx_{\al_0}^0+ (-\ep
r_2\be_1-\be_0)(N^0)^Tdx_{\be_0}^0= i\frac{\ep_1(\ep
r_2\be_0-\be_1)}{\sqrt{H_0}}d\al_1
-i\frac{\ep_1(r_1\al_0-\al_1)}{\sqrt{H_0}}d\be_1$, or using
$\frac{\pa}{\pa\al_0}=\frac{\mathbf{C}_0}{\la_0}\frac{\pa}{\pa
u}-\frac{\mathbf{S}_0}{\mu_0}\frac{\pa}{\pa v},\
\frac{\pa}{\pa\be_0}=-\ep\frac{\mathbf{S}_0}{\la_0}\frac{\pa}{\pa
u}+\frac{\mathbf{C}_0}{\mu_0}\frac{\pa}{\pa v},\
-\ep(r_1\al_1-\al_0)d\al_0+(-\ep r_2\be_1-\be_0)d\be_0=
\ep(r_1\al_0-\al_1)d\al_1+(\ep r_2\be_0-\be_1)d\be_1$:
$\\\mu_0[(r_1\al_1-\al_0)\mathbf{C_0}-\ep(-\ep
r_2\be_1-\be_0)\mathbf{S_0}]=
\ep_1[(\ep r_2\be_0-\be_1)\al_{1u}-(r_1\al_0-\al_1)\be_{1u}],\\
\la_0[(-\ep
r_2\be_1-\be_0)\mathbf{C_0}-(r_1\al_1-\al_0)\mathbf{S_0}]=
-\ep_1[(r_1\al_0-\al_1)\be_{1v}-(\ep r_2\be_0-\be_1)\al_{1v}],\\
\la_0[(r_1\al_1-\al_0)\mathbf{C_0}-\ep(-\ep
r_2\be_1-\be_0)\mathbf{S_0}]=
-[(r_1\al_0-\al_1)\al_{1u}+\ep(\ep r_2\be_0-\be_1)\be_{1u}],\\
\mu_0[(-\ep
r_2\be_1-\be_0)\mathbf{C_0}-(r_1\al_1-\al_0)\mathbf{S_0}]= [(\ep
r_2\be_0-\be_1)\be_{1v}+\ep(r_1\al_0-\al_1)\al_{1v}]$.

If these equations are $I-IV$ and using $(\ep
r_2\be_0-\be_1)^2+\ep(r_1\al_0-\al_1)^2=-zH_0>0$, then by
considering $I^2+\ep III^2,\ IV^2+\ep II^2,\ I\cdot II-III\cdot
IV$ we obtain $\al_{1u}^2+\ep\be_{1u}^2>0,\
\be_{1v}^2+\ep\al_{1v}^2>0,\\
\al_{1u}\al_{1v}+\ep\be_{1u}\be_{1v}=0$, so
$\begin{bmatrix}\al_{1u}&\al_{1v}\\\be_{1u}&\be_{1v}\end{bmatrix}=
\begin{bmatrix}\la_1\mathbf{C}_1&-\ep\mu_1\mathbf{S}_1\\
\la_1\mathbf{S}_1&\mu_1\mathbf{C}_1\end{bmatrix}$. Because of the
symmetry $(0,\ep)\leftrightarrow(1,-\ep)$ and using $I,\ II$ we
obtain that the second fundamental form of $x^1$ is
$-i\frac{\la_1\mu_1(du^2+dv^2)}{\sqrt{H_1}}$, so $(u,v)$ is also
isothermic-conjugate on $x^1$ and the B transformation changes the
orientation of $(u,v)$ (at the level of the analytic computations
this can be accounted by putting $\ep_0:=-\ep_1$ with $N^1$ also
changing sign); from $II^2-\ep III^2$ we obtain
$\mu_1^2-\ep\la_1^2=H_1$ and we have complete symmetry
$(0,\ep)\leftrightarrow(1,-\ep)$ also at the level of the
isothermic-conjugate system $(u,v)$. Thus
$\\\al_1=r_1\al_0-\ep_1'\sqrt{-z}(\mathbf{C}_1\la_0+\ep\ep_1\mathbf{S}_1\mu_0),\
\be_1=\ep
r_2\be_0-\ep_1'\sqrt{-z}(\mathbf{S}_1\la_0-\ep_1\mathbf{C}_1\mu_0),\
(\ep_0,\ep'_0):=-(\ep_1,\ep'_1),\\
(0,\ep)\leftrightarrow(1,-\ep)$.

Finally differentiating these with respect to $u,v$ we obtain that
the B transformation (\ref{eq:backesgh}) between solutions
$\om_0,\ \ep_1\theta_1$ and respectively $\theta_0,\ \ep_1\om_1$
of (\ref{eq:esg}) and (\ref{eq:esgh}) has as influence the
algebraic transformations of solutions of (\ref{eq:linsyste})
respectively (\ref{eq:linsystie}) as follows:

\subsection{Real deformations of (the imaginary region of) the
real elliptic paraboloid} \noindent

\noindent

\subsection{Imaginary deformations of (the imaginary
region of) the real elliptic paraboloid} \noindent

\noindent

\subsection{Real deformations of (the imaginary region of) the real hyperboloid with one
sheet} \noindent

\noindent

Consider the general confocal real hyperboloids with one sheet in
an isothermic-conjugate parametrization invariant under the Ivory
affinity between confocal quadrics
$$x_z=x_z(\al,\be):=[\sqrt{a_1-z}\cos\be\sec\al\ \ \sqrt{a_2-z}\sin\be\sec\al\ \
\sqrt{-a_3+z}\tan\al]^T,\
a_1>a_2>z,0>a_3,$$$$\al\in(-\frac{\pi}{2},\frac{\pi}{2}),\be\in\mathbb{R}$$
with linear element, second fundamental form and Christoffel
symbols of $x_0$
$$|dx_0|^2=a_1(d(\cos\be\sec\al))^2
+a_2(d(\sin\be\sec\al))^2-a_3(d\tan\al)^2;$$
$$N_0^Td^2x_0=\frac{d\al^2- d\be^2}{\cos\al\sqrt{H}},\
H:=a_1^{-1}\cos^2\be+a_2^{-1}\sin^2\be-a_3^{-1}\sin^2\al;$$
$$-\Ga_{22}^1=\Ga_{11}^1-2\tan\al=(\log\sqrt{H})_{\al},\
-\Ga_{11}^2=\Ga_{22}^2=(\log\sqrt{H})_{\be},\ \Ga_{12}^2=\tan\al,\
\Ga_{12}^1=0.$$ We also have a distinguished tangent vector field
$\mathcal{V}_0:=((\log\sqrt{H})_{\al}+\tan\al)x_{0\al}-(\log\sqrt{H})_{\be}x_{0\be}$;
it has the properties
\begin{eqnarray}\label{eq:v0}
dx_{0\al}=\mathcal{V}_0d\al-N_0dN_0^Tx_{0\al}-\tan\al dx_0,\nonumber\\
dx_{0\be}=-\mathcal{V}_0d\be-N_0dN_0^Tx_{0\be}+\tan\al(x_{0\be}d\al
+x_{0\al}d\be).
\end{eqnarray}
Given a real deformation $x\subset\mathbb{R}^3$ of a real region
$\subset x_0$ there exists a conjugate system $(u,v)$ common to
both $x_0$ and $x$ (this is true for any two surfaces in a
point-wise correspondence). Denote with $\bar\cdot$ the quantities
of interest in the GCMP equations (namely the Christoffel symbols
and the second fundamental form) of $x_0$ referred to the $(u,v)$
coordinates and similarly with $\tilde\cdot$ those of $x$. We have
$\al_u\al_v-\be_u\be_v=0$ and from the Gau\ss\ equation
$(\al_u^2-\be_u^2)(\al_v^2-\be_v^2)<0$; assume
$\al_u^2-\be_u^2>0$. With
$\la:=\mathrm{sgn}(\al_u)\sqrt{\al_u^2-\be_u^2},\
\mu:=\mathrm{sgn}(\be_v)\sqrt{\be_v^2-\al_v^2}$ we have $\bar
h_{11}=\frac{\la^2}{\cos\al\sqrt{H}},\ \bar h_{12}=0,\ \bar
h_{22}=-\frac{\mu^2}{\cos\al\sqrt{H}}$. From the general formula
for the change of Christoffel symbols $\frac{\pa u^l}{\pa\ti
u^c}\ti\Ga_{ab}^c=\frac{\pa^2u^l}{\pa\ti u^a\pa\ti u^b}+\frac{\pa
u^j}{\pa\ti u^a}\frac{\pa u^k}{\pa\ti u^b}\Ga_{jk}^l$ we are
interested only in $\bar\Ga_{12}^1,\ \bar\Ga_{12}^2$:
$\al_{uv}+2\tan\al\al_u\al_v=\bar\Ga_{12}^1\al_u+\bar\Ga_{12}^2\al_v,\
\be_{uv}+\tan\al(\al_u\be_v+\al_v\be_u)=\bar\Ga_{12}^1\be_u+\bar\Ga_{12}^2\be_v$,
so $\bar\Ga_{12}^1=(\log\frac{\la}{\cos\al})_v,\
\bar\Ga_{12}^2=(\log\frac{\mu}{\cos\al})_u$. From the CMP
equations of $x_0,\ x$ we have
$$(\bar h_{11})_v=\bar\Ga_{12}^1\bar h_{11}-\bar\Ga_{11}^2\bar
h_{22},\ (\ti h_{11})_v=\bar\Ga_{12}^1\ti h_{11}-\bar\Ga_{11}^2\ti
h_{22},$$
$$(\bar h_{22})_u=\bar\Ga_{12}^2\bar h_{22}-\bar\Ga_{22}^1\bar
h_{11},\ (\ti h_{22})_u=\bar\Ga_{12}^2\ti h_{22}-\bar\Ga_{22}^1\ti
h_{11}.$$ Keeping account of the Gau\ss\ equation $\bar h_{11}\bar
h_{22}=\ti h_{11}\ti h_{22}$ one can multiply the first equations
respectively with $\bar h_{11},\ \ti h_{22}$ (and the second
equations respectively with $\bar h_{22},\ \ti h_{11}$), subtract
them and get rid respectively of the $\bar\Ga_{11}^2,\
\bar\Ga_{22}^1$ terms: $(\log(\bar h_{11}^2-\ti
h_{11}^2))_v=2(\log\frac{\la}{\cos\al})_v,\ (\log(\bar
h_{22}^2-\ti h_{22}^2))_u=2(\log\frac{\mu}{\cos\al})_u$. Thus
$\bar h_{11}^2-\ti h_{11}^2=\phi(u)\frac{\la^2}{\cos^2\al},\ \bar
h_{22}^2-\ti h_{22}^2=\varphi(v)\frac{\mu^2}{\cos^2\al}$; after a
change of the $u$ and $v$ variables one can absorb $\phi(u),\
\varphi(v)$ up to opposite signs $\ep_1:=\pm 1$: $\ti
h_{11}^2=\bar
h_{11}^2-\ep_1\frac{\la^2}{\cos^2\al}=\frac{\la^2}{\cos^2\al}(\frac{\la^2}{H}-\ep_1),\
\ti h_{22}^2=\bar
h_{22}^2+\ep_1\frac{\mu^2}{\cos^2\al}=\frac{\mu^2}{\cos^2\al}(\frac{\mu^2}{H}+\ep_1)$
(again from the Gau\ss\ equation we have $\frac{|\bar
h_{11}|}{|\ti h_{11}|}=\frac{|\ti h_{22}|}{|\bar h_{22}|}$). Now
from the Gau\ss\ equation $\ti h_{11}\ti h_{22}=\bar h_{11}\bar
h_{22}$ we get $H=-\ep_1\mu^2+\ep_1\la^2$; we can choose
$\ep_1:=-1$. Thus the second fundamental form of $x$ is
$\frac{\la\mu(du^2-dv^2)}{\cos\al\sqrt{H}}$ and we are led to
consider the hyperbolic angle $\theta$ between the conjugate
systems $(\al,\be)$ and $(u,v)$, that is
$\begin{bmatrix}\al_u&\al_v\\\be_u&\be_v\end{bmatrix}=
\begin{bmatrix}\la\mathbf{C}&\mu\mathbf{S}\\\la\mathbf{S}&\mu\mathbf{C}\end{bmatrix},\
\mathbf{C}:=\cosh\theta,\ \mathbf{S}:=\sinh\theta$ (note that by
doing this the sign of $\theta$ will be decided by that of
$\be_u$). Imposing the compatibility conditions
$(\la\mathbf{C})_v=(\mu\mathbf{S})_u,\
(\la\mathbf{S})_v=(\mu\mathbf{C})_u$ we get $\mu_u=\la\theta_v,\
\la_v=\mu\theta_u$; differentiating $H=\mu^2-\la^2$ with respect
to $u$, respectively $v$ we are led to consider the differential
system in $\al,\ \be,\ \la,\ \mu$:
\begin{eqnarray}\label{eq:linsysth}
d\begin{bmatrix}\al\\\be\\\la\\\mu\end{bmatrix}=\begin{bmatrix}
\mathbf{C}\la du+\mathbf{S}\mu dv\\\mathbf{S}\la du+\mathbf{C}\mu
dv\\
[-\frac{1}{2}H_{\al}\mathbf{C}-\frac{1}{2}H_{\be}\mathbf{S}
+\mu\theta_v]du+\mu\theta_udv\\
\la\theta_vdu+[\frac{1}{2}H_{\al}\mathbf{S}+\frac{1}{2}H_{\be}\mathbf{C}+
\la\theta_u]dv\end{bmatrix},\ \mu^2-\la^2=H
\end{eqnarray}
with the modified hyperbolic sinh-Gordon equation
\begin{eqnarray}\label{eq:mhypSGh}
\theta_{vv}-\theta_{uu}=\frac{1}{2}(H_{\al\al}+
H_{\be\be})\mathbf{C}\mathbf{S}
\end{eqnarray}
as the compatibility condition. As opposed to the case of
paraboloids, where $H$ was quadratic, (\ref{eq:mhypSGh}) has the
inconvenience of depending on $\al,\be$ and the dependence of
$\al,\be$ on $(u,v)$ being undetermined, so further manipulation
is required to obtain an equation depending on $\theta$ only (it
will be a third or fourth order differential equation), but it is
not our interest to do that now; note also that its B
transformation will reveal itself later.

Note that if we assume that the common conjugate system on $x_0,\
x$ is isothermic-conjugate on $x$ (Darboux), then from the Gau\ss\
equations we obtain immediately that the second fundamental form
of $x$ is $\frac{\la\mu(du^2-dv^2)}{\cos\al\sqrt{H}}$; everything
else except $\mu_v,\ \la_u,\ H=\mu^2-\la^2$ follows immediately as
previously. The remaining needed information follows immediately
from the CMP equations of $x$: we have
$\bar\Ga_{12}^1=(\log\frac{\la}{\cos\al})_v,\
\bar\Ga_{22}^1=\frac{\mu}{\la}\theta_v-\frac{\mu^2}{\la^2}(\log\sqrt{H})_u,\
\bar\Ga_{11}^2=\frac{\la}{\mu}\theta_u-\frac{\la^2}{\mu^2}(\log\sqrt{H})_v,\
\bar\Ga_{12}^2=(\log\frac{\mu}{\cos\al})_u$ and
$(\log\frac{\la\mu}{\cos\al\sqrt{H}})_v=\bar\Ga_{12}^1+\bar\Ga_{11}^2,\
(\log\frac{\la\mu}{\cos\al\sqrt{H}})_u=\bar\Ga_{12}^2+\bar\Ga_{22}^1$
respectively become
$\mu_v=\frac{\mu^2-\la^2}{\mu}(\log\sqrt{H})_v+\la\theta_u,\
\la_u=-\frac{\mu^2-\la^2}{\la}(\log\sqrt{H})_u+\mu\theta_v$, so
$d\log\frac{\mu^2-\la^2}{H}=0$; by a same homothety in the $(u,v)$
variables we and a choice of sign we can assume
$\frac{\mu^2-\la^2}{H}=1$.

Next we derive the algebraic computations of the {\it tangency
configuration} (TC) and of the B transformation.

Note
\begin{eqnarray}\label{eq:basic}
-x_{z\al}x_{z\al}^T+x_{z\be}x_{z\be}^T+\sec^2\al
x_zx_z^T=\sec^2\al\mathrm{diag}[a_1-z\ \ a_2-z\ \ a_3-z].
\end{eqnarray}
Consider two points $x_0^0,\ x_0^1\in x_0$; with
$V_0^1:=x_z^1-x_0^0,\ \hat N_0:=-2\pa_z|_{z=0}x_z$ the relation
\begin{eqnarray}\label{eq:baseb}
(-x_{z\al_1}^1(x_{z\al_1}^1)^T+x_{z\be_1}^1(x_{z\be_1}^1)^T)\hat
N_0^0=-\sec^2\al_1[x_z^1(V_0^1)^T\hat N_0^0+(V_0^1+z\hat N_0^0)]
\end{eqnarray}
follows immediately from (\ref{eq:basic}).

Consider now the symmetric TC $(V_0^1)^T\hat N_0^0=(V_1^0)^T\hat
N_0^1=0$. Multiplying (\ref{eq:baseb}) on the left respectively
with $(\hat N_0^0)^T,\ (V_0^1\times N_0^0)^T$ we get two algebraic
consequences of the TC:
\begin{eqnarray}\label{eq:actc}
-[(x_{z\al_1}^1)^TN_0^0]^2+[(x_{z\be_1}^1)^TN_0^0]^2=-z\sec^2\al_1=
\frac{|dx_z^1|^2-|dx_0^1|^2}{-d\al_1^2+d\be_1^2},\nonumber\\
(x_{z\al_1}^1-x_{z\be_1}^1)^T(I_3-2N_0^0(N_0^0)^T)[(x_{z\al_1}^1+x_{z\be_1}^1)\times
V_0^1]=0.
\end{eqnarray}
Note that $x_{z\al}\pm x_{z\be}$ are rulings on $x_z$ and thus
their length is preserved under the Ivory affinity.

With $(R_0^1,t_0^1)$ being the {\it rigid motion provided by the
Ivory affinity} (RMPIA) such that $\\(R_0^1,t_0^1)[x_0^0\ \ x_z^1\
\ x_{0\al_0}^0-x_{0\be_0}^0\ \ x_{z\al_1}^1-x_{z\be_1}^1]=[x_z^0\
\ x_0^1\ \ x_{z\al_0}^0-x_{z\be_0}^0\ \
x_{0\al_1}^1-x_{0\be_1}^1]$ we have
$x_{0\al_1}^1+x_{0\be_1}^1=R_0^1(I_3-2N_0^0(N_0^0)^T)(x_{z\al_1}^1+x_{z\be_1}^1)$
since by changing the ruling family on $x_z^1$ the action of the
RMPIA on $T_{x_0^0}x_0$ does not change.

Multiplying (\ref{eq:baseb}) on the left with $R_0^1$ and using
(\ref{eq:actc}) we obtain
\begin{eqnarray}\label{eq:ivor}
x_z^1=x_0^0+\cos^2\al_0(-x_{0\al_0}^0(x_{z\al_0}^0)^T+x_{0\be_0}^0(x_{z\be_0}^0)^T)\hat
N_0^1,\ 0\leftrightarrow 1.
\end{eqnarray}
Differentiating the relation $(V_1^0)^T\hat N_0^1=0$ we obtain
\begin{eqnarray}\label{eq:difftc}
(dx_z^0)^T\hat N_0^1=-(x_z^0)^Td\hat N_0^1=-(\hat N_0^0)^Tdx_z^1.
\end{eqnarray}
Now Bianchi's main theorem on the deformation of quadrics (the
existence and inversion of the B\"{a}cklund transformation and the
ACPIA) states (in our case) that the TC coupled with
\begin{eqnarray}\label{eq:back}
x^1=x^0+\cos^2\al_0(-x_{\al_0}^0(x_{z\al_0}^0)^T+x_{\be_0}^0(x_{z\be_0}^0)^T)\hat
N_0^1,\ |dx^1|^2=|dx_0^1|^2
\end{eqnarray}
is a differential system in involution (completely integrable)
given the {\it seed} deformation $x^0\subset\mathbb{R}^3$ of
$x_0^0$ (that is $|dx^0|^2=|dx_0^0|^2$), that the $1$-dimensional
family of solutions ({\it leaves}) $x^1$ is given by the
integration of a Ricatti equation, that $x^0$ and $x^1$ are the
focal surfaces of a Weingarten congruence (congruence of lines on
whose two focal surfaces the asymptotic directions correspond;
since conjugate directions are harmonically conjugate to the
asymptotic ones all conjugate systems correspond in this case) and
that we have the symmetry $0\leftrightarrow 1$.

If $(R_0,t_0)\subset\mathbf{O}_3(\mathbb{R})\ltimes\mathbb{R}^3$
is the rolling of $x_0^0$ on $x^0$ (that is
$(R_0,t_0)(x_0^0,dx_0^0):=(R_0x_0^0+t_0,R_0dx_0^0)=(x^0,dx^0)),\
N^0:=R_0N_0^0$, then $R_0^{-1}dR_0N_0^0=R_0^{-1}dN^0-dN_0^0,\
R_0^{-1}dx^1=dx_z^1+R_0^{-1}dR_0V_0^1=dx_z^1+(V_0^1)^TR_0^{-1}dR_0N_0^0N_0^0
=(I_3-N_0^0(N_0^0)^T)dx_z^1-N_0^0(dN^0)^T(x^1-x^0)$ and
$|dx_0^1|^2=|dx^1|^2$ becomes
$[(dN^0)^T(x^1-x^0)]^2=|dx_0^1|^2-|dx_z^1|^2+[(N_0^0)^Tdx_z^1]^2=
[(N_0^0)^T(x_{z\be_1}^1d\al_1+x_{z\al_1}^1d\be_1)]^2$, so
\begin{eqnarray}\label{eq:factor}
(dN^0)^T(x^1-x^0)=\ep_1(N_0^0)^T(x_{z\be_1}^1d\al_1+x_{z\al_1}^1d\be_1),\
\ep_1:=\pm 1
\end{eqnarray}
(each choice of $\ep_1$ corresponds to a ruling family). Using
$\begin{bmatrix}\frac{\pa}{\pa\al_0}\\\frac{\pa}{\pa\be_0}\end{bmatrix}=
\begin{bmatrix}\frac{\mathbf{C}_0}{\la_0}&-\frac{\mathbf{S}_0}{\mu_0}\\
-\frac{\mathbf{S}_0}{\la_0}&\frac{\mathbf{C}_0}{\mu_0}\end{bmatrix}
\begin{bmatrix}\frac{\pa}{\pa
u}\\\frac{\pa}{\pa v}\end{bmatrix}$, (\ref{eq:difftc}),
(\ref{eq:back}) and (\ref{eq:factor}) we get:
$\\\mu_0[\mathbf{C_0}(x_{z\al_0}^0)^T+\mathbf{S_0}(x_{z\be_0}^0)^T]\hat
N_0^1=
\ep_1[\al_{1u}(x_{z\be_1}^1)^T+\be_{1u}(x_{z\al_1}^1)^T]\hat N_0^0,\\
\la_0[\mathbf{C_0}(x_{z\be_0}^0)^T+\mathbf{S_0}(x_{z\al_0}^0)^T]\hat
N_0^1=
\ep_1[\be_{1v}(x_{z\al_1}^1)^T+\al_{1v}(x_{z\be_1}^1)^T]\hat N_0^0,\\
\la_0[\mathbf{C_0}(x_{z\al_0}^0)^T+\mathbf{S_0}(x_{z\be_0}^0)^T]\hat
N_0^1=
-[\al_{1u}(x_{z\al_1}^1)^T+\be_{1u}(x_{z\be_1}^1)^T]\hat N_0^0,\\
\mu_0[\mathbf{C_0}(x_{z\be_0}^0)^T+\mathbf{S_0}(x_{z\al_0}^0)^T]\hat
N_0^1=-[\be_{1v}(x_{z\be_1}^1)^T+\al_{1v}(x_{z\al_1}^1)^T]\hat
N_0^0$.

Consider the case $z<0$. If these equations are $I-IV$ and using
(\ref{eq:actc}) (here we use $z<0$), then by considering
$I^2-III^2,\ IV^2-II^2,\ I\cdot II-III\cdot IV$ we obtain
$\al_{1u}^2-\be_{1u}^2>0,\ \be_{1v}^2-\al_{1v}^2>0,\
\al_{1u}\al_{1v}-\be_{1u}\be_{1v}=0$, so
$\begin{bmatrix}\al_{1u}&\al_{1v}\\\be_{1u}&\be_{1v}\end{bmatrix}=
\begin{bmatrix}\la_1\mathbf{C}_1&\mu_1\mathbf{S}_1\\
\la_1\mathbf{S}_1&\mu_1\mathbf{C}_1\end{bmatrix}$. Because of the
symmetry $0\leftrightarrow 1$ and using $I,\ II$ we obtain that
the second fundamental form of $x^1$ is
$\frac{\la_1\mu_1(du^2-dv^2)}{\cos\al_1\sqrt{H_1}}$, so $(u,v)$ is
also isothermic-conjugate on $x^1$ and the B transformation
preserves the orientation of $(u,v)$ (for $z>0$ it changes it);
from $II^2-III^2$ we obtain $\mu_1^2-\la_1^2=H_1$ and we have
complete symmetry $0\leftrightarrow 1$ also at the level of the
isothermic-conjugate system $(u,v)$. Thus $\\(x_{z\al_0}^0)^T\hat
N_0^1=\ep_1'\sqrt{-z}\sec\al_0\sec\al_1(\mathbf{C}_0\la_1+\ep_1\mathbf{S}_0\mu_1),\\
(x_{z\be_0}^0)^T\hat
N_0^1=-\ep_1'\sqrt{-z}\sec\al_0\sec\al_1(\mathbf{S}_0\la_1+\ep_1\mathbf{C}_0\mu_1),\
\ep_0':=-\ep_1',\ 0\leftrightarrow 1$.

With $r_j:=\sqrt{1-za_j^{-1}}$ we have the TC
$r_1\cos\be_0\cos\be_1+r_2\sin\be_0\sin\be_1=\cos\al_0\cos\al_1+r_3\sin\al_0\sin\al_1$
and the above become $\\\sin\al_0\cos\al_1-r_3\cos\al_0\sin\al_1=
\ep_1'\sqrt{-z}(\mathbf{C}_0\la_1+\ep_1\mathbf{S}_0\mu_1),\\
-r_1\sin\be_0\cos\be_1+r_2\cos\be_0\sin\be_1=
-\ep_1'\sqrt{-z}(\mathbf{S}_0\la_1+\ep_1\mathbf{C}_0\mu_1),\
\ep_0':=-\ep_1',\ 0\leftrightarrow 1$, or
$\\\begin{bmatrix}\sin\al_0\cos\al_1\\\cos\al_0\sin\al_1\end{bmatrix}
=\frac{-a_3}{\ep_1'\sqrt{-z}}\begin{bmatrix}1&r_3\\r_3&1\end{bmatrix}
\begin{bmatrix}\mathbf{C}_0\la_1+\ep_1\mathbf{S}_0\mu_1\\
-(\mathbf{C}_1\la_0+\ep_1\mathbf{S}_1\mu_0)\end{bmatrix},\\
\begin{bmatrix}\sin\be_0\cos\be_1\\\cos\be_0\sin\be_1\end{bmatrix}=
-\frac{1}{\ep_1'\sqrt{-z}(a_2^{-1}-a_1^{-1})}\begin{bmatrix}r_1&r_2\\r_2&r_1\end{bmatrix}
\begin{bmatrix}\mathbf{S}_0\la_1+\ep_1\mathbf{C}_0\mu_1\\
-(\mathbf{S}_1\la_0+\ep_1\mathbf{C}_1\mu_0)\end{bmatrix}$.

We have
$\\(\mathbf{C}_0\la_1+\ep_1\mathbf{S}_0\mu_1)_u=(\mathbf{S}_0\la_1+\ep_1\mathbf{C}_0\mu_1)
(\theta_{0u}+\ep_1\theta_{1v})-\mathbf{C}_0(\frac{1}{2}H_{1\al_1}\mathbf{C}_1
+\frac{1}{2}H_{1\be_1}\mathbf{S}_1),\\
(\mathbf{C}_0\la_1+\ep_1\mathbf{S}_0\mu_1)_v=(\mathbf{S}_0\la_1+\ep_1\mathbf{C}_0\mu_1)
(\theta_{0v}+\ep_1\theta_{1u})+\ep_1\mathbf{S}_0(\frac{1}{2}H_{1\al_1}\mathbf{S}_1
+\frac{1}{2}H_{1\be_1}\mathbf{C}_1),\\
(\mathbf{S}_0\la_1+\ep_1\mathbf{C}_0\mu_1)_u=(\mathbf{C}_0\la_1+\ep_1\mathbf{S}_0\mu_1)
(\theta_{0u}+\ep_1\theta_{1v})-\mathbf{S}_0(\frac{1}{2}H_{1\al_1}\mathbf{C}_1
+\frac{1}{2}H_{1\be_1}\mathbf{S}_1),\\
(\mathbf{S}_0\la_1+\ep_1\mathbf{C}_0\mu_1)_v=(\mathbf{C}_0\la_1+\ep_1\mathbf{S}_0\mu_1)
(\theta_{0v}+\ep_1\theta_{1u})+\ep_1\mathbf{C}_0(\frac{1}{2}H_{1\al_1}\mathbf{S}_1
+\frac{1}{2}H_{1\be_1}\mathbf{C}_1)$ and the B transformation of
(\ref{eq:mhypSGh}) in conjunction with solutions of
(\ref{eq:linsysth}) reveals itself:
$\\\theta_{0u}+\ep_1\theta_{1v}=\frac{\ep_1'}{\sqrt{-z}}[
(\sin\al_0\sin\al_1+r_3\cos\al_0\cos\al_1)\mathbf{S}_0\mathbf{C}_1-
(r_1\sin\be_0\sin\be_1+r_2\cos\be_0\cos\be_1)\mathbf{C}_0\mathbf{S}_1],\\
\theta_{0v}+\ep_1\theta_{1u}=\frac{-\ep_1\ep_1'}{\sqrt{-z}}[
(\sin\al_0\sin\al_1+r_3\cos\al_0\cos\al_1)\mathbf{C}_0\mathbf{S}_1-
(r_1\sin\be_0\sin\be_1+r_2\cos\be_0\cos\be_1)\mathbf{S}_0\mathbf{C}_1]$

\subsection{Calapso's B\"{a}cklund transformation for real
quadrics of revolution} \noindent

\noindent

\subsection{Darboux's integral formula for deformations
of the real paraboloid of revolution} \noindent

\noindent

\section{The solitons of quadrics}

In analogy to the situation for the link between the solitons of
the sine-Gordon equation and the solitons of the pseudo-sphere
(when the $0$-soliton is the axis of the tractrix) we are
interested in finding degenerate deformations of quadrics (that is
the seed collapses to a curve or point) as $0$-solitons and then
in finding explicit formulae of their B transforms.

For real deformations of the real hyperbolic paraboloid from the
differential system (\ref{eq:linsyst}) we have $\la=0,\
\al=\al(v),\ \be=\be(v),\ \mu=\mu(v),\ \theta=\theta(v),\
\al'=\mu\mathbf{S},\ \be'=\mu\mathbf{C},\
-\mathbf{C}\frac{\al}{a_1}+\mathbf{S}\frac{\be}{a_2}+\mu\theta'=0,\
\mu'=\mathbf{S}\frac{\al}{a_1}-\mathbf{C}\frac{\be}{a_2},\
\mu^2=\frac{\al^2}{a_1}-\frac{\be^2}{a_2}+1$ and $\theta$ will
satisfy the (hyperbolic) pendulum equation
$\theta''=\mathbf{S}\mathbf{C}$.

We have $(\theta')^2=\frac{\mathbf{C}^2+\mathbf{S}^2+c}{2}$, so
the solution $\theta$ is given in terms of elliptic functions
$v=\ep\int\frac{\sqrt{2}d\theta}{\sqrt{\mathbf{C}^2+\mathbf{S}^2+c}},\\
\ep=\pm 1$. Then we take $\al$ solution of the second order ODE
$\al''-2(\log\mathbf{S})'\al'+\frac{\al}{a_1}=0$ and
$\be:=-a_2(\frac{\theta'}{\mathbf{S}^2}\al'-\frac{\mathbf{C}}{\mathbf{S}}\frac{\al}{a_1}),\
\mu:=\frac{1}{\theta'}(\mathbf{C}\frac{\al}{a_1}-\mathbf{S}\frac{\be}{a_2})$
and we have $\al'=\mu\mathbf{S},\ \be'=\mu\mathbf{C},\
\mu'=\mathbf{S}\frac{\al}{a_1}-\mathbf{C}\frac{\be}{a_2},\
\mu^2-\frac{\al^2}{a_1}+\frac{\be^2}{a_2}=\mathrm{ct}$; this last
constant can be normalized to $1$ by a choice of constant in the
initial value of $\al$. Note the ODE of $\al$ can be brought to
the form $(\mathbf{C}^2+\mathbf{S}^2+c)\frac{d^2\al}{d\theta^2}-
2\frac{\mathbf{C}}{\mathbf{S}}(\mathbf{C}^2+c)\frac{d\al}{d\theta}+2\frac{\al}{a_1}=0$
and the prime integral above mentioned to the form
$\frac{\mathbf{C}^2+\mathbf{S}^2+c}{2\mathbf{S}^2}(\frac{d\al}{d\theta})^2+
\frac{a_1}{1-a_1}(\frac{\mathbf{C}^2+\mathbf{S}^2+c}{2\mathbf{S}^2}\frac{d\al}{d\theta}
-\frac{\mathbf{C}}{\mathbf{S}}\frac{\al}{a_1})^2-\frac{\al^2}{a_1}=1$.
The homogeneous part can be integrated by quadrature
($\frac{d\log\al}{d\theta}$ will depend algebraically on
$\mathbf{C},\mathbf{S}$) and then by the standard variation of
parameters argument one can solve this prime integral for
$\al=\al(\theta)$.

Thus all quantities of interest can be found by explicit formulae;
the only remaining question is if one can find explicit formulae
for the B transforms of the solutions of the hyperbolic pendulum
equation.

Again a change from the $v$ variable to the $\theta_0$ variable is
in order ($\theta_1=\theta_1(u,\theta_0)$):

$$\theta_{1\theta_0}=\ep\frac{\sqrt{2}}{\sqrt{\mathbf{C}_0^2+\mathbf{S}_0^2+c}}
(\frac{\si_1+\si_1^{-1}}{2}\mathbf{S}_1\mathbf{C}_0
+\frac{\si_1-\si_1^{-1}}{2}\mathbf{C}_1\mathbf{S}_0),$$
$$\theta_{1u}=\ep\frac{\sqrt{\mathbf{C}_0^2+\mathbf{S}_0^2+c}}{\sqrt{2}}
+\frac{\si_1-\si_1^{-1}}{2}\mathbf{S}_1\mathbf{C}_0
+\frac{\si_1+\si_1^{-1}}{2}\mathbf{C}_1\mathbf{S}_0.$$ The last
equation is separable (we consider $\theta_0=\mathrm{ct}$); by
quadrature one can find the solution depending on a constant of
$u$ (function of $\theta_0$); in turn by replacing the result in
the first equation one can find the function of $\theta_0$ up to a
constant.

Thus the $0$-solitons will depend on two constants and each
iteration of the B transformation will introduce two constants.

For $c=\pm 1$ the elliptic function will degenerate to hyperbolic
trigonometric ones and $\theta_0$ will turn out to be the
$1$-solitons of the hyperbolic sinh-Gordon equation with $\si=1$;
then one can apply directly the BPT to find $\theta_1$.

Also since we already know the $1$-solitons of the hyperbolic
sinh-Gordon equation from Peterson's deformations of quadrics, a
space realization of solitons is possible in this particular case.

\end{document}